\documentclass[a4paper,10pt]{article}
\usepackage[top=1 in,bottom=1in, left=1.5 in, right=1.5 in]{geometry}

\usepackage[utf8]{inputenc}
\usepackage{graphicx}
\usepackage{amsmath}
\usepackage{authblk}
\usepackage{hyperref}       
\usepackage{url}            
\usepackage{booktabs}       
\usepackage{amsfonts}       
\usepackage[ugly]{nicefrac}       
\usepackage{microtype}      
\usepackage{lipsum}
\usepackage{csquotes}

\usepackage{array}
\usepackage{multicol}	
\usepackage{float}
\usepackage{amsthm}
\usepackage{amsmath} 
\usepackage{amssymb}
\usepackage{bbm}
\usepackage{epstopdf}
\usepackage{color,fancybox,graphicx}
\usepackage{url}
\usepackage{psfrag}
\usepackage{upgreek}
\usepackage{pstricks}
\usepackage{graphicx,psfrag}
\usepackage[dvips]{epsfig}
\usepackage{array}
\usepackage{enumitem}
\usepackage{float}
\usepackage{latexsym}

\newcommand{\x}{\textbf{x}} 
\newcommand{\y}{\textbf{y}} 
\newcommand{\xx}{\textbf{X}} 
\newcommand{\s}{S}
\newcommand{\Sn}{S_n}
\newcommand{\Rn}{R_n}
\newcommand{\E}{\mathbb{X}}
\newcommand{\set}{\mathbb{X}}

\newcommand{\argmin}{\operatornamewithlimits{argmin}}
\newcommand{\voo}{\boldsymbol{\xi}_{n}} 
\newcommand{\doe}{\textbf{D}_n}  
\newcommand{\gvoo}{\textbf{g}_n}
\newcommand{\loiX}{\textbf{P}_{\xx}}
\newcommand{\bornemarkov}{\delta_n}
\newcommand{\borneconvexe}{\gamma_n}
\newcommand{\JR}{J_{\Rn}}
\newcommand{\Ji}{J_{1,n}}
\newcommand{\Jii}{J_{2,n}}
\newcommand{\Jiii}{J_{3,n}}
\newcommand{\Jiv}{J_{4,n}}
\newcommand{\JS}{J_{\Sn}}
\newcommand{\Jk}{J_{k,n}}

\newcommand{\tribu}{\mathcal{F}}

\theoremstyle{theorem}
\newtheorem{rem}{Remark}[section]
\newtheorem{Pte}{Property}[section]
\newtheorem{defi}{Definition}[section]

\newtheorem{pro}{Proposition}[section]

\begin{document}

\title{Estimating a probability of failure with the convex order in computer experiments}
\author[1]{Lucie~Bernard\thanks{Email: lucie.bernard-3@etu.univ-tours.fr}}
\author[2]{Philippe Leduc\thanks{Email: philippe.leduc@st.com}}
\date{June 26, 2019\\\small{}}
\affil[1]{\small{Institut Denis Poisson, Universit\'e de Tours, Tours, France}}
\affil[2]{\small{ STMicroelectronics, Tours, France}}
\maketitle

\begin{abstract}
This paper deals with the estimation of a failure probability of an industrial product. To be more specific, it is defined as the probability that the output of a physical model, with random input variables, exceeds a threshold. The model corresponds with an expensive to evaluate black-box function, so that classical Monte Carlo simulation methods cannot be applied. Bayesian principles of the Kriging method are then used to design an estimator of the failure probability. From a numerical point of view, the practical use of this estimator is restricted. An alternative estimator is proposed, which is equivalent in term of bias. The main result of this paper concerns the existence of a convex order inequality between these two estimators. This inequality allows to compare their efficiency and to quantify the uncertainty on the results that these estimators provide. A sequential procedure for the construction of a design of computer experiments, based on the principle of the Stepwise Uncertainty Reduction strategies, also results of the convex order inequality. The interest of this approach is highlighted through the study of a real case from the company STMicroelectronics.
\end{abstract}

\section{Introduction}

To evaluate the profitability of a production before the launch of its manufacturing process, most industrial companies use computer simulation in order to measure the impact of the variability inherent to it. It consists in virtually testing possible configurations of the fluctuating parameters of the given industrial product, and predicting its performance by studying the simulation results in light of specifications. In this paper, the product under study is defined by a number of~$d\in\mathbb{N}^*$ fluctuating parameters. Generally called factors, these parameters vary in a set $\set\subset\mathbb{R}^d$, assumed to be known. Moreover, each simulation is supposed to require the execution of a complex and expensive calculation code, so that we consider the model: \begin{align*}y=g(\x),\end{align*} where the vector input $\x\in\set$ is a set of factors that represents some experimental conditions and 
$g:\E\to\mathbb{R}$ is a black-box function that is costly to evaluate. It is assumed deterministic, which means that evaluations of this function with the same input always return the same output. The scalar output $y$ is a numerical value, generally called response, that measures the product performance in the configuration describes by~$\x$. 
Besides, a probability distribution~$\loiX$ is given on  $\set$ and reflects the variability of the factors. We assume that we know how to simulate according to this distribution, so that a sequence $(\xx_i)_{1\leq i\leq N}$, where $N\in\mathbb{N}^*$, of random variables i.i.d. according to~$\loiX$ is available. 
Here, a failure consists in the output of a computer experiment exceeding a threshold $T\in\mathbb{R}$, so that 
 we are interested in estimating the probability $p$ defined by:\begin{align*}
p = \mathbb{P}\left(g(\xx)>T\right)=\displaystyle\int_{\E}\mathbbm{1}_{g(\x)>T}\textbf{P}_\xx(d\x),
\end{align*} 
 usually called {failure probability} in industrial fields. 
The main difficulty is obviously that the evaluation budget for the function $g$ needed to estimate $p$ is severely limited. That is why 
a naive Monte Carlo estimation method, which consists in simulating an i.i.d. sample $(\xx_i)_{1\leq i\leq N}$ w.r.t. $\textbf{P}_\xx$ and setting: \begin{align*} \widehat{p}_N=\dfrac{1}{N}\sum_{i=1}^N\mathbbm{1}_{g(\xx_i)>T},
\end{align*} 
is clearly unreasonable. Indeed, recall that the accuracy of such an estimate is usually given by the relative standard deviation of $\widehat{p}_N$, which is equal to $\sqrt{(1-p)/Np}$. This means that, for instance, one needs to select a sample of size $N$ of order $10^5$ in order to obtain a relative precision equals to $10\%$ for the estimation of a probability of $10^{-3}$. Since any evaluation of $g$ is very expensive, such a computational burden is unrealistic. Thus, we assume in all this paper that only a budget of $n\ll N$ evaluations of $g$ is available in order to estimate $p$.  
This implies, on the one hand, setting a design of experiments $(\x_i)_{1\leq i \leq n}\in\set^n$ where to perform the evaluations and, on the other hand, building a reliable estimator despite this limited number of observations. In an industrial context, it means that we aim at choosing the configurations of the product parameters to be explored first by numerical simulation in order to obtain relevant information on its performance. Moreover, we are interested in providing a measure of uncertainty about the result obtained, which is essential when information is scarce and not systematically provided by methods already existing in the literature to answer a similar problem.\medskip

The function $g$ being unknown, we focus in this paper on 
a modelling method called \enquote{Kriging}. Initially developed for geostatistics (see, e.g., \cite{Diggle07}, \cite{Matheron63}, \cite{Stein99}), this approach gained later considerations in the design and analysis of computers experiments (see, e.g., \cite{Rasmussen06}, \cite{Sacks89_DA} and \cite{TheDesign}). The applications to failure probability estimation and, more generally, structural reliability problems are rather recent (see, e.g., \cite{EGR}, \cite{Dubourg13} and \cite{Krig_Struc_Realiability_Pbs}). In \cite{BMA14}, Bayesian principles of Kriging are used to design an estimator of the probability of failure $p$, defined as a random variable $\Sn$ taking values in $[0,1]$. From a Bayesian point of view, the limitations of this estimator are as follows: its exact distribution is out of reach, its empirical distribution is difficult to calculate and its caracteristics (Bayes risk, variance, quantiles, \ldots) cannot be estimated at a lower cost.  
In \cite{Oger}, an alternative estimator is proposed, which is equivalent in terms of bias. 
Denoted by $\Rn$, we aim at providing methods to learn about the distribution of $\Sn$ by considering $\Rn$. Our main contribution is about the existence of a convex order between these two estimators, meaning that for all convex function $\varphi$, we have: \begin{align}\label{jdheys}
\mathbb{E}\left[\varphi(\Rn)\right]\leq\mathbb{E}\left[\varphi(\Sn)\right].
\end{align} We show that Inequality (\ref{jdheys}) allows for comparing their efficiency and can be exploited to achieve our objectives. Since we consider that the mean value of $\Sn$ provides an estimation of $p$, we build credible intervals in order to quantify the uncertainty on this measurement. For a fixed value of $\alpha\in[0.1]$, an interval $[a_\alpha, b_\alpha]$ satisfying:\begin{align}\label{sstszjeus}
\mathbb{P}(a_\alpha\leq \Sn\leq b_\alpha)\geq 1-\alpha,
\end{align}is a $(1-\alpha)$-credible  interval for the distribution of $\Sn$. 
The credible intervals we propose here are derived from Inequality (\ref{jdheys}) and give better results than the one based on the Markov's inequality proposed in \cite{BMA14}. Among all the methods already available in the literature which also use the Kriging method to answer a problem of estimating a probability of failure, this information is rarely available. The second concerns the choice of the experimental design $(\x_i)_{1\leq i \leq n}$. Based on the principle of the Stepwise Uncertainty Reduction (SUR) strategies presented in \cite{Li12}, we use our result on the convex order to implement a sequential procedure. By using a criterion based on the variance of $\Rn$, the points to be added to the design of experiments are iteratively selected and the model is adjusted according to these new observations. The interest of this approach is particularly highlighted through the study of a real case from the company STMicroelectronics. As a reference, a Monte Carlo analysis based on 1000 numerical simulation results have been performed and we will see that the estimations obtained with a much smaller number of simulations are compatible with those obtained with 1000.\medskip

In Section \ref{Section2}, we first briefly recalls the Kriging method in a Bayesian framework. Then we define the random variable $S_n$, and we recall the advantages and limitations of this approach encountered in previous work on the subject.
In Section \ref{Section3}, we present the random variable $\Rn$ and, in section  \ref{Section4}, we show that there is a convex order between these two random variables.  In Section \ref{Section5}, several applications of this result are given and we are particularly interested in the construction of credibility intervals (\ref{sstszjeus}). In addition, we describe in Section \ref{Section6} the principles of {SUR} strategies for the construction of the experimental design and introduce our alternative sampling criterion. In Section \ref{Section7}, we propose some illustrations for a one-dimensional artificial example and an industrial case. Finally, a conclusion and perspectives are given in Section \ref{Section8}. Proofs are detailed in Section \ref{Section9}.

\section{A first estimator}\label{Section2}



Let us assume that function $g$ is observed at points of a fixed design of experiments $\doe = (\x_i)_{1\leq i\leq n}\in\E^n$ and   denote by $\gvoo=(g(\x_i))_{1\leq i\leq n}$ the vector of observations. We precise that we do not address here the problem of the choice of the most efficient design to perform evaluations of $g$. We refer instead to \cite{EGR}, \cite{Dette}, \cite{Sacks89_D} and the sequential strategy presented in section \ref{Section5}.

		\subsection{Gaussian process models}
		
Bayesian principles of the Kriging method are basically the following. First, the function~$g$ is assumed to be a sample path of a random process $\xi$ indexed by $\set$, which is generally chosen Gaussian (prior model). This is indeed one of the rare model for which analytical calculations can be carried out and which is well suited to a number of observed physical phenomena. Then, the distribution of~$\xi$ is conditioned by the fact that a realization of the random vector $\xi_{\doe}  = (\xi(\x_i))_{1\leq i \leq n}$ is observed and equal to~$g_{\doe}$. It turns out that the process is still Gaussian (posterior model). All its sample paths interpolate the observed points (for more information, see, e.g., \cite{Rasmussen06}, \cite{TheDesign} and \cite{Stein99}). 
To simplify, we note~$\xi_n$ the random process with the same finite-dimensional laws as the posterior process $\xi$ conditioned on observations, which implies the following equality in law: \begin{align}\label{qtgehtudte}
\mathcal{L}(\xi(\x)\mid \voo = \gvoo)=\mathcal{L}(\xi_n(\x)), \quad\forall\x\in\set.
\end{align}For all $\x\in\set$, we denote by $m_n(\x)  =  \mathbb{E}[\xi(\x) \mid \voo=\gvoo]=\mathbb{E}[\xi_n(\x)]$ the Kriging mean and $\sigma^2_n(\x)=\mbox{Var}[\xi(\x)\mid \voo =\gvoo]=\mbox{Var}[\xi_n(\x)]$ the Kriging variance, so that $\xi_n(\x)\thicksim\mathcal{N}(m_n(\x), \sigma_n^2(\x))$, $\forall\x\in\E$. For all $i=1,\ldots,n$, they satisfy $m_n(\x_i)=g(\x_i)$ and $\sigma_n^2(\x_i)=0$, which means that all sample paths of $\xi_n$ are indeed functions that interpolate points of coordinates $(\doe, g_{\doe})$. Explicit expressions of the Kriging mean and Kriging variance can be found in many books and publications (see, e.g., \cite{Sacks89_DA} and references mentioned above). They are not useful in the rest of this paper though, so it is not considered necessary to recall them here. However, we point out that these quantities are easily computable in practice, provided that the covariance matrix $(\mbox{Cov}(\xi(\x_i), \xi(\x_j)))_{1\leq i, j\leq n}$ is positive-definite. Thus, it depends on the crucial choice of the covariance function that characterizes the dependency between the points of the design $\doe$. See, e.g., \cite{Abrahamsen97} for general properties and \cite{Rasmussen06} for an exhaustive list of covariance functions, but also standard methods for selecting one based on the observations (the so-called Leave-One-Out cross-validation method, for instance). Note that the results presented further are completely independent of the choice of the covariance function. Besides, although it is more convenient to choose a Gaussian random process, our results are actually valid for any choice of a prior distribution, provided that one can easily compute the following probability of membership to the failure set:
\begin{align}\label{pn}
p_n(\x) 
= \mathbb{P}\big(\xi(\x)>T\mid \voo=\gvoo\big) = \mathbb{P}\big(\xi_n(\x)>T\big),\quad\forall\x\in\E,
\end{align}which satisfies in the Gaussian case:\begin{align*}p_n(\x) = \Phi\left(\frac{m_n(\x) - T}{\sigma_n(\x)}\right),\quad\forall\x\in\E,\end{align*}
where $\Phi$ is the cumulative distribution function  of the standard Gaussian distribution. We precise in advance that one-dimensional examples of Gaussian process modeling are given Figure \ref{g_Ex1}. The dashed black line is the function $m_n$ and gray areas are $95\%$-confidence intervals written as follows:\begin{align}\label{kdolemmmmm}
\left[m_n(\x) - 1.96\sigma_n(\x)\quad;\quad m_n(\x) + 1.96\sigma_n(\x)\right], \quad\forall\x\in\E.
\end{align}
	
		\subsection{Bayesian approach}\label{First_Estimate}
			
In order to build an estimator of $p$, Bayesian principles of Kriging can be used as follows. Under the assumption that $g$ is a trajectory of a Gaussian process $\xi$, the probability $p$ is a realization of the random variable $\s\in[0,1]$ defined by:
\begin{equation}\label{S}
\s = \mathbb{P}\left(\xi(\xx) > T | \xi\right) = \displaystyle\int_{\mathbb{X}}\mathbbm{1}_{\xi(\x)> T} \textbf{P}_{\xx}(d\x).
\end{equation}
Intuitively, the probability that a trajectory of $\xi$ exceeds the threshold $T$ is a realization of $\s$. As stated in \cite{BMA14}, it is straightforward to show that the posterior distribution of $\s$ knowing that $\voo = \gvoo$ is the same as the random variable~$\Sn\in[0,1]$ defined by:\begin{equation}\label{R1}
\Sn = \mathbb{P}\left(\xi_n(\xx) > T | \xi_n\right) = \displaystyle\int_{\mathbb{X}}\mathbbm{1}_{\xi_n(\x)> T} \textbf{P}_{\xx}(d\x).
\end{equation}If the process $\xi_n$ is sufficiently informed so that sample paths accurately reproduce the behaviour of $g$ in the neighborhood of $T$,  it is then reasonable to assume that any realization of $\Sn$ is relatively close to~$p$. Consequently, the first estimator we consider in this paper is the random variable $\Sn$ and the quantities of interest are, in particular, its mean value, its variance and its quantiles. Thus, let us denote by $\mu_n$ the mean value of~$\Sn$, i.e. the posterior mean value of $\s$ knowing $\voo = \gvoo$. It is well-known that, in the Bayesian framework, this quantity refers to the Bayes estimator of $\s$ with respect to the quadratic loss function (see, e.g., \cite{Livre_Choix_Bayes}). Here, $\mu_n$ has a very simple analytic expression and is written: \begin{align}\label{Mean_R1}
\mu_n = \mathbb{E}\left[\Sn\right] &  = \displaystyle\int_{\mathbb{X}}\mathbb{E} [\mathbbm{1}_{\xi_n(\x)> T}] \textbf{P}_{\xx}(d\x) 
  =  \displaystyle\int_{\mathbb{X}} p_n(\x) \textbf{P}_{\xx}(d\x) = \mathbb{E}[p_n(\xx)].
\end{align}
The function $p_n$ given in Equation (\ref{pn}) is rather inexpensive to evaluate, so that we can apply a naive Monte Carlo method to estimate $\mu_n$. This does not require additional calls to the function $g$. Indeed, by taking $N\in\mathbb{N}^*$ large enough, it is sufficient to generate an i.i.d. sequence of random variables $(\xx_i)_{1\leq i \leq N}$ with respect to $\loiX$ and consider the following approximation:\begin{align}\label{lodkdh}
\mu_n \approx \frac{1}{N}\sum_{i=1}^Np_n(\xx_i).
\end{align}
The numerical value, returned by (\ref{lodkdh}) in practice, can be used to estimate the probability of failure $p$. 

	\subsection{Limits}\label{qhqye}

To evaluate the efficiency of the estimator $\Sn$, it is desirable to calculate its variance. In the Bayesian framework, this quantity refers to the posterior risk with respect to the quadratic loss function (see, e.g., \cite{Livre_Choix_Bayes}). One can easily check that it is equal to:
\begin{align}\label{Var_Sn}
\mbox{Var}\left[\Sn\right]= \displaystyle\int_{\mathbb{X}^2}\mathbb{P}\big(\xi_n(\x_1)> T, \xi_n(\x_2)> T\big) \textbf{P}_{\xx}(d\x_1)\textbf{P}_{\xx}(d\x_2) - \mu_n^2.
\end{align} Note that, more generally, the raw moment of order $m\in\mathbb{N}^*$ satisfies: 
\begin{align}\label{Second_moment_R1}
\mathbb{E}[\Sn^m]  & = \mathbb{E}\left[\displaystyle\int_{\E^m}\bigg(\prod_{i=1}^m\mathbbm{1}_{\xi_n(\x_i)> T}\bigg) \textbf{P}_{\xx}(d\x_1)\ldots\textbf{P}_{\xx}(d\x_m)\right]\nonumber\\
				   & = \displaystyle\int_{\mathbb{X}^m}\mathbb{P}\bigg(\bigcap_{i=1}^m \big\{\xi_n(\x_i)> T\big\}\bigg) \textbf{P}_{\xx}(d\x_1)\ldots\textbf{P}_{\xx}(d\x_m).
\end{align}In practice, these measures can be estimated by using a naive Monte Carlo method, an estimator of the right term in (\ref{Var_Sn}) is already given in Equation (\ref{lodkdh}). Nevertheless, the computation time requested by this approach can be high since it involves the computation of a joint probability. Besides, according to our knowledge, to get the distribution of $\Sn$ is unachievable and it seems unreasonable to learn about it by making statistical inference from the analysis of realizations of $\Sn$. Indeed, generate a realization of this random variable involves, on the one hand, simulating a trajectory of the random process $\xi_n$ and, on the other hand, performing an integration with respect to $\textbf{P}_{\xx}$ (through a naive Monte Carlo method, for instance). A method based on the discretization of the set $\E$ for the simulation of trajectories conditionally on observations is proposed in \cite{BMA14}. However, the authors acknowledge that this approach could lead to burdensome calculations, since it requires the inversion of a covariance matrix whose size depends on the number of points selected for discretization. Yet, the latter grows exponentially with the dimension of $\E$. Note that authors in \cite{BMA14} also propose to upper-bound the quantiles of $\Sn$ with an approach based on the Markov's inequality, and another one based on the principle of the Importance Sampling method, which is a variance reduction technique for the Monte Carlo method. The first one will be detailed in Section \ref{Section5}. 


\section{An alternative estimator}\label{Section3}	

From a Bayesian point of view, it is natural to consider the random variable $\Sn$ introduced in the last section, despite its distribution is out of reach. To learn about it and, consequently, justify the interest of the Bayesian approach for estimating the probability of failure $p$, we use an alternative estimation method proposed in \cite{Oger}. 
By taking a random variable $U$ uniform on $[0,1]$, it involves the alternative random variable $\Rn\in[0,1]$ defined~by: \begin{equation}\label{R2}
\Rn = \mathbb{P}(p_n(\xx) > U | U) = \displaystyle\int_{\mathbb{X}}\mathbbm{1}_{p_n(\x) > U}\textbf{P}_{\xx}(d\x),
\end{equation} 
where $p_n$ is defined in Equation (\ref{pn}). It is easy to verify that $\Sn$ and $\Rn$ are equivalent in terms of bias:\begin{align*}
\mathbb{E}[\Rn] =  \displaystyle\int_{\mathbb{X}}\mathbb{E}[\mathbbm{1}_{p_n(\x) > U}]\textbf{P}_{\xx}(d\x) = \int_\E p_n(\x)\loiX(d\x) = \mu_n.
\end{align*}
Note that, in \cite{Oger}, the authors prove that the uniform law on $[0,1]$ is the only one for $U$ such that, for any law $\loiX$ and for any distribution choice for the $\xi$ process, we have:\begin{align*}\mathbb{E}[\Rn] = \mathbb{E}[\Sn] = \mu_n.
\end{align*}Except for this result of equality of the means, the relationship between $\Sn$ and $\Rn$ is not further studied in \cite{Oger}, the authors having being focused on intrinsically justify the choice of $\Rn$ for the estimation of $p$, in case of a deterministic function $g$. The second estimator we consider in this paper is then the random variable $\Rn$ and, at this stage, its advantages are the  following. Firstly, the law of $\Rn$ does not directly depend on the law of $\xi_n$.  Indeed, the information provided by the fact that $\voo = \gvoo$ is taken into account through the function $p_n$ and ensures that only the marginal laws of $\xi_n$ are involved in the modelling. Secondly, the function $p_n$ has the advantage of being relatively inexpensive to evaluate and the random variable $U$ is unidimensional, regardless of the $d$ dimension of $\E$. Finally, unlike $\Sn$, 
we can approximately simulate according to the distribution of $\Rn$ at a lower cost. Indeed, if we generate an i.i.d. sample $(U_j)_{1\leq j\leq M}$ w.r.t. the standard uniform law and another one $(\xx_i)_{1\leq i\leq N}$ w.r.t. $\loiX$, then we can consider that:\begin{align*}
\int_{\E}\mathbbm{1}_{p_n(\x) > U_j}\textbf{P}_{\xx}(d\x)\approx\dfrac{1}{N}\sum_{i=1}^N\mathbbm{1}_{p_n(\xx_i)>Uj},\quad\forall j=1, \ldots, M.
\end{align*}

\section{A convex order inequality}\label{Section4}

The theory of stochastic orders aims at providing tools to compare the distributions of random variables to assist in decision making. An overview of the many existing orders can be found in \cite{StochasticOrders07} and a basic introduction of the most popular stochastic orders is given in \cite{MulSto2002}. The convex order is one of them and is usually used to compare the variability of random variables which have the same mean value. As a reminder, the definition of the convex order is the following (see, e.g., \cite{StochasticOrders07}, Chapter 3): 
\begin{defi}\label{def_cx_order}
The random variable X is said to be smaller than the random variable Y in the convex order, denoted $X\leq_{cx}Y$, if for all convex function $\varphi:\mathbb{R}\to\mathbb{R}$, provided that the expectations exist, we have:
$$\mathbb{E}[\varphi(X)]\leq\mathbb{E}[\varphi(Y)].$$
\end{defi}  
\noindent
Recall that $\Sn$ and $\Rn$ are equivalent in terms of bias. Our main result concerns the existence of such a relationship between these random variables, meaning that it is possible to compare their efficiency.
\begin{pro}\label{prop_cx_order}
The random variables $\Sn$ and $\Rn$ defined in (\ref{R1}) and (\ref{R2}) satisfy: \begin{align*}
\Sn\leq_{cx}\Rn.
\end{align*}
\end{pro}\noindent
According to Definition \ref{def_cx_order}, we immediately notice that Proposition \ref{prop_cx_order} implies that:\begin{align}\label{lqueje}
\mathbb{E}[\Sn^m] \leq \mathbb{E}[\Rn^m], ~\forall m\in\mathbb{N}^*,\quad\mbox{and}\quad\mbox{Var}[\Sn]\leq\mbox{Var}[\Rn].
\end{align}Recall that the variance and the $m$-th moment of $\Sn$ are respectively given in Equations (\ref{Var_Sn}) and (\ref{Second_moment_R1}). Then, Inequalities (\ref{lqueje}) are in fact immediate, since we have:\begin{align}\label{Var_R2}
\mathbb{E}[\Rn^m] & = \displaystyle\int_0^1 \left(\displaystyle\int_\mathbbm{X}\mathbbm{1}_{p_n(\x) > u} \textbf{P}_{\xx}(d\x)\right)^mdu \nonumber\\
				   & = \displaystyle\int_{\mathbb{X}^m}\left(\displaystyle\int_0^1\mathbbm{1}_{p_n(\x_1) > u} \ldots\mathbbm{1}_{p_n(\x_m) > u} du\right)\textbf{P}_{\xx}(d\x_1)\ldots\textbf{P}_{\xx}(d\x_m) \nonumber\\
				   & = \displaystyle\int_{\mathbbm{X}^2}\min\big(p_n(\x_1),\ldots ,p_n(\x_m)\big)\textbf{P}_{\xx}(d\x_1)\ldots\textbf{P}_{\xx}(d\x_m),
\end{align} and $\mathbb{P}\left(A\cap B\right) \leq\min\left(\mathbb{P}(A), \mathbb{P}(B)\right)$. 
The Inequalities  (\ref{lqueje}) show that 
 $\Sn$ is a better estimator than $\Rn$ in terms of precision. 
Nevertheless, we give priority to $\Rn$ because, as mentioned above, it has the same mean value $\mu_n$ and is proved to be of better practical use. Besides, Inequalities (\ref{lqueje}) also ensure that we can estimate bounds of the raw moments of $\Sn$. Indeed, unlike
suggested by Equation (\ref{Var_R2}), the calculation cost to estimate the raw moments of $\Rn$ does not increase with the order $m$. Since we can simulate a sequence of random variables approximately i.i.d.~according to the distribution of $\Rn$ (see Section \ref{Section3}), it is sufficient to calculate the empirical moments. Concerning the variance of $\Rn$, more details  will be given in Section \ref{Section6}. Before that, let us precise that many properties related to the convex order are verified by $\Sn$ and $\Rn$ (see Chapter 3 of \cite{StochasticOrders07} for an overview). We focus, in the next section, on those that provide information on the distribution of $\Sn$ and approximations of its quantiles. 

\section{Bounding quantiles using convex order}\label{Section5}


	\subsection{Improved bounds}

Recall that the quantile function $F_X^{-1}$ of a random variable $X$ is defined for all $\alpha\in [0,1]$ by:\begin{align}\label{zyehtgdtre}
F_X^{-1}(\alpha) = \inf\lbrace t\in\mathbb{R} : \mathbb{P}( X \leq t ) \geq \alpha\rbrace,
\end{align}
where the quantity $F_X^{-1}(\alpha)$ is called $\alpha$-quantile of $X$. It is not always easy to estimate a quantile and, generally, the proposed approach is non-parametric: an i.i.d.~sample is simulated according to the law of $X$ and the empirical estimator is considered (see, e.g., \cite{Livre_stat}). 
There are many applications in finance that refer to the quantile of a distribution, used as a measure of  risk called Value-at-Risk (see, e.g., \cite{Livre_VaR}, \cite{Consistency_risks} et \cite{Estimation_risks}). In the following proposition, we show that the quantile function of~$\Rn$ has in fact an analytical expression, so that it can be estimated using a naive Monte Carlo method based on simulations w.r.t. $\loiX$ (see Section \ref{First_Estimate}). 		
\begin{pro}\label{Prop_CDF_Quantile_function_Rn}
Let $\Rn$ be the random variable defined in Equation (\ref{R2}). Then the quantile function $F_{\Rn}^{-1}$ of $\Rn$ satisfies:\emph{\begin{align*} F_{\Rn}^{-1}(\alpha) = \int_\E\mathbbm{1}_{p_n(\x)> 1- \alpha}\textbf{P}_{\xx}(d\x), \quad\forall\alpha\in [0,1].\end{align*}}
\end{pro}
\noindent
Now, we aim at providing an approximation of the quantile function of $\Sn$. For this purpose, authors in \cite{BMA14} suggest to use the inequalities of Markov and Chebyshev's inequality. 
Here, they lead to the following:\begin{align}\label{Borne_Markov_Barbillon}
F_{\Sn}^{-1}(\alpha)\leq \dfrac{\mu_n}{1-\alpha}
\quad\mbox{and}\quad F_{\Sn}^{-1}(\alpha)\leq\mu_n+\sqrt{\dfrac{\mbox{Var}[\Sn]}{1-\alpha}}, \quad\forall\alpha\in (0,1).
\end{align}In Section \ref{First_Estimate}, we have explained how to estimate $\mu_n$ at low cost, so that the first bound in (\ref{Borne_Markov_Barbillon}) is easy to estimate. On the other hand -- and this is underlined in \cite{BMA14} -- the Monte Carlo method for the estimation of the variance of $\Sn$ requires relatively high computation times (this point has already been discuss in Section \ref{qhqye}). A simple way to get around the problem is to bound from above the variance of $\Sn$ in (\ref{Borne_Markov_Barbillon}) by the variance of $\Rn$:\begin{align}\label{njnjn}
F_{\Sn}^{-1}(\alpha)\leq\mu_n+\sqrt{\dfrac{\mbox{Var}[\Rn]}{1-\alpha}}.
\end{align} The variance of $\Rn$ is cheap to estimate (we recall that an easy-to-implement method will be proposed in Section \ref{Implementation}). Nevertheless, we acknowledge that the upper-bound (\ref{njnjn}) based on the variance of $\Rn$ provides a less effective control than the one based on the variance of $\Sn$. That is why, in order to achieve a compromise between precision and computational complexity, the following bounds are proposed. 
	
\begin{pro}\label{Prop_CredibleInt_Markov}
For all  $\alpha\in (0,1)$, the quantile $F_{\Sn}^{-1}(\alpha)$ satisfies: 
\begin{align}\label{Inequality_CredibleInt_Markov}
\bornemarkov^-(\alpha)\leq\borneconvexe^-(\alpha)\leq F_{\Sn}^{-1}(\alpha)\leq\borneconvexe^+(\alpha)\leq\bornemarkov^+(\alpha),
\end{align}where:\begin{align*}
\bornemarkov^-(\alpha) = \frac{\mu_n+\alpha-1}{\alpha}\quad\mbox{and}\quad\bornemarkov^+(\alpha) =\dfrac{\mu_n}{1-\alpha},
\end{align*}and\begin{align*}
\borneconvexe^-(\alpha) =\dfrac{1}{\alpha}\displaystyle\int_0^\alpha F_{\Rn}^{-1}(t)dt \quad\mbox{and}\quad\borneconvexe^+(\alpha)=\dfrac{1}{1 - \alpha}\displaystyle\int_\alpha^1 F_{\Rn}^{-1}(t)dt.
\end{align*}
\end{pro}
\noindent
To our best knowledge, the intermediate bounds $\borneconvexe^-(\alpha)$ and $\borneconvexe^+(\alpha)$ given in Proposition \ref{Prop_CredibleInt_Markov} are new to address this problematic.
 In the proof, we show that they are obtained by using the convex order inequality between $\Sn$ et $\Rn$, so that will call them \enquote{convex order bounds}. The bounds $\bornemarkov^-(\alpha)$ and $\bornemarkov^+(\alpha)$ refer to the Markov's inequality (\ref{Borne_Markov_Barbillon}) and therefore we will call them  \enquote{Markov's bounds}.  
In practice, we should take $\max(0, \bornemarkov^-(\alpha))$ and $\min(1, \bornemarkov^+(\alpha))$  since they are not systematically informative. Indeed, 
\begin{align*}
\mu_n\leq 1-\alpha\Rightarrow \bornemarkov^-(\alpha)\leq 0\leq F_{\Sn}^{-1}(\alpha),
\end{align*}and \begin{align*}
{\mu_n}\geq 1-\alpha\Rightarrow F_{\Sn}^{-1}(\alpha)\leq 1\leq \bornemarkov^+(\alpha).
\end{align*}
 
This restriction is not verified for the convex order bounds, because they necessarily take values in $(0,1)$. Besides, we show in the following proposition that these latter can be written as an integrand w.r.t. the law $\loiX$, meaning that a naive Monte Carlo method can again be apply.

\begin{pro}\label{Exp_bounds}
For all  $\alpha\in (0,1)$, let $\borneconvexe^-(\alpha)$ and $\borneconvexe^+(\alpha)$ be the bounds given in Proposition \ref{Prop_CredibleInt_Markov}. Then, we have:
\emph{\begin{align*}\borneconvexe^-(\alpha) = 1-\displaystyle\int_{\E}\min\left(1, \dfrac{1-p_n(\x)}{\alpha}\right)\textbf{P}_{\xx}(d\x) ,\end{align*}}
and\emph{\begin{align*}\borneconvexe^+(\alpha) = \displaystyle\int_{\E}\min\left(1, \dfrac{p_n(\x)}{1-\alpha}\right)\textbf{P}_{\xx}(d\x).\end{align*}}
\end{pro}According to Propositions \ref{Prop_CredibleInt_Markov} and \ref{Exp_bounds}, there is therefore no practical difficulty for approximating the quantiles of~$\Sn$. 

\begin{rem}
In the financial and actuarial science literature, the quantity $\frac{1}{1-\alpha}\int_{\alpha}^1F_{X}^{-1}(t)dt$ is a risk measure for $X$, usually called {Conditional Value-at-Risk} at level $\alpha\in[0,1)$  (see \cite{MulSto2002} and \cite{MRA}, Section 7.1.2, and \cite{ESB}). For more information on its properties and the related estimation methods, see also \cite{CTE_applicationsActurial}. 
\end{rem}

	\subsection{Credible intervals}

For a fixed value $\alpha\in(0,1)$, our goal is here to determine an interval $[a_\alpha,b_\alpha]\subseteq[0,1]$ satisfying:\begin{align}\label{bbreqtgr}
\mathbb{P}(a_\alpha\leq \Sn\leq b_\alpha)\geq 1-\alpha.\end{align}In the Bayesian framework, such an interval is called  $(1-\alpha)$-credible interval (or credible region). A formal definition can be found in \cite{Livre_Choix_Bayes}, for instance. We obviously want this interval to be easy to estimate, but it also needs to be narrow in order to provide relevant information. For any value of $\alpha\in[0,1]$, there are indeed an infinite number of intervals~$[a_\alpha,b_\alpha]$ satisfying (\ref{bbreqtgr}). According to Proposition \ref{Prop_CredibleInt_Markov}, the convex order's bounds~$\borneconvexe^-(\alpha)$ and $\borneconvexe^+(\alpha)$ introduced in the previous section satisfy:
\begin{align*}
\borneconvexe^-(\alpha)\leq F_{\Sn}^{-1}(\alpha)\leq \borneconvexe^+(\alpha).
\end{align*}Therefore, considering also Equation \ref{zyehtgdtre}, we have for all $\beta\in (0,1)$:\begin{align*}
\mathbb{P}\big(\Sn\leq\borneconvexe^-(\alpha\beta)\big)\leq {\alpha}{\beta}, 
\end{align*}and\begin{align*}
\mathbb{P}\big(\Sn\leq\borneconvexe^+(1-\alpha(1-\beta))\big)\geq 1 - {\alpha}{(1-\beta)}, 
\end{align*}In other words,  \begin{align*}
\mathbb{P}\big(\borneconvexe^-(\alpha\beta)\leq\Sn\leq \borneconvexe^+(1-\alpha(1-\beta))\big)\geq 1-\alpha, 
\end{align*}
The following proposal summarizes this and uses the expressions of $\borneconvexe^-(\alpha)$ and $\borneconvexe^+(\alpha)$ given in Proposition \ref{Exp_bounds}.
\begin{pro}\label{dsjsdjzDJsef}
Let $\alpha\in (0,1)$ be fixed. For all $\beta\in (0,1)$, the interval $I^{cx}_n(\alpha, \beta)$ defined by:\emph{\begin{align}\label{Intervalle_cred_cx}
I^{cx}_n(\alpha, \beta) = \left[1-\int_{\mathbb{X}}\min\left(1,  \dfrac{1-p_n(\x)}{\alpha\beta}\right)\textbf{P}_{\xx}(d\x)\quad,\quad \int_{\mathbb{X}}\min\left(1, \dfrac{p_n(\x)}{\alpha(1-\beta)}\right)\textbf{P}_{\xx}(d\x)\right],
\end{align}}satisfies:\begin{align*}
\mathbb{P}\big(\Sn\in I^{cx}_n(\alpha, \beta)\big)\geq 1-\alpha, \quad\forall\beta\in (0,1).
\end{align*}
\end{pro}


For a fixed value of $\alpha\in(0,1)$, the length of the interval given in Equation (\ref{Intervalle_cred_cx}) depends on the value of the parameter~$\beta\in(0,1)$. We then suggest to use in practice a classical optimization algorithm to determine the value of~$\beta$ which minimizes it. Otherwise, we can simply take $\beta=\frac{1}{2}$, for instance. Note that the interval $I^{Markov}_n(\alpha,\beta)$ defined for all $\beta\in(0,1)$ by:\begin{align}\label{Intervalle_cred_markov}
I^{Markov}_n(\alpha,\beta) = \left[\frac{\mu_n + \alpha\beta - 1}{\alpha\beta}\quad,\quad\dfrac{\mu_n}{\alpha(1-\beta)}\right],
\end{align} also satisfies (\ref{bbreqtgr}). According to Proposition \ref{Prop_CredibleInt_Markov}, it is nevertheless  wider than the one given in Equation~(\ref{Intervalle_cred_cx}).

\section{Application of the convex order to sequential design of computer experiments}\label{Section6}

When the credible interval specified in Proposition  \ref{dsjsdjzDJsef} is too large to consider (\ref{lodkdh}) as a reliable estimation of the probability of failure  $p$, we conclude that the Kriging model is not sufficiently informed. To improve the prediction quality, it is necessary to add information to the model, that is to provide new observations. To this end, we propose a sequential procedure for selecting the design of experiments. It is based on the principle of Stepwise Uncertainty Reduction (SUR) strategies, which have been formalized in \cite{Li12} for the Kriging framework. For a better understanding, we simplify here the formalism and adapt it to our study.  For a more general presentation of {SUR} strategies, we can of course consult \cite{Li12}, but also the thesis \cite{Chevalier13_Thesis} and references therein. In addition, we precise that theoretical results justifying the performance of these methods have recently been proposed in \cite{bect_martingale}. 

	\subsection{SUR strategies}
	
	For all $n\in\mathbb{N}^*$, we suppose that the function $g$ has been already observed at points of a design of experiments $\doe = (\x_i)_{1\leq i \leq n}$ and we aim at finding the point $\x^*_{n+1}\in\set\setminus\doe$ where perform the next evaluation of $g$. For this purpose, we provide the probability space $(\Omega, \tribu, \mathbb{P})$ with 
the $\sigma$-algebra $\tribu_n$ generated by the random variables $(\xi(\x_i))_{1\leq i \leq n}$ for any set of points $(\x_i)_{1\leq i \leq n}$. 
In these conditions, $\xi_n$ is the random process with the same finite-dimensional distributions as $\xi$ knowing $\tribu_n$, the random variable $\Sn$ is distributed according to the law of $\s$ knowing $\mathcal{F}_n$ and $\Rn$ is the random variable defined by $
\Rn = \int_{\set}\mathbbm{1}_{p_n(\x)>U}\loiX(d\x)$, where $p_n(\x) =\mathbb{P}(\xi(\x)>T\mid \mathcal{F}_n)=\mathbb{P}(\xi_n(\x)>T)$, $\forall\x\in\set$. In addition, we adopt the following notation: $\mathbb{E}_n [\cdot] = \mathbb{E}[\cdot\vert\tribu_n],\forall n \in\mathbb{N}^*$.
		\subsubsection{Minimization of a criterion}


When $n$ observations are available, the basic principle of SUR strategies is to define an uncertainty measure $H_{n+1}$ on a quantity of interest, depending on the objective to achieve, and to select the point $\x_{n+1}^*$ that decreases the most this uncertainty. Here, we work in the particular setting where the quantity of interest is $S$ and its uncertainty measure is its posterior variance: \begin{align}\label{H_n}
H_{n+1} = \mbox{Var}[S\mid \mathcal{F}_{n+1}] = \mbox{Var}[S_{n+1}].
\end{align} It is desirable for this variance to be as small as possible, because it ensures the relevance of the quantity (\ref{lodkdh}) to estimate $p$ when $n+1$ observations will be available. 
Nevertheless, at this stage, the variance (\ref{H_n}) is a random variable that depends on $\xi(\x_{n+1})$, for which only the information is the following: \begin{align*}\mathcal{L}(\xi(\x_{n+1})\mid \tribu_n) = \mathcal{N}(m_{n}(\x_{n+1}), \sigma_n^2(\x_{n+1})).\end{align*} Therefore, it seems reasonable to select the point $\x_{n+1}^*$ minimizing the variance of $S_{n+1}$ in expectation:
\begin{align*}\x_{n+1}^* = \argmin_{\x_{n+1}\in\mathbb{X}} \mathbb{E}_n\big[\mbox{Var}[S_{n+1}]\mid\xx_{n+1} = \x_{n+1}\big].\end{align*}This leads to the definition of the following sampling criterion, which is here noted as~$\JS$, and whose minimization leads to identify where to perform the next evaluation of  $g$:\begin{equation*}\JS(\x_{n+1}) = \mathbb{E}_n\big[\mbox{Var}[S_{n+1}]\mid\xx_{n+1} = \x_{n+1}\big].
\end{equation*}

		\subsubsection{Alternative criteria}

In practice, it is difficult to evaluate criterion $\JS$ without requiring high computation times (see Section \ref{qhqye}). In \cite{Li12}, the authors   acknowledge the limitations of this approach and propose alternative criteria. For all $\x\in \E$, let us first define: \begin{align*}
\tau_n(\x) = \min(p_n(\x), 1-p_n(\x)) \quad\mbox{and}\quad\nu_n(\x) = p_n(\x)(1-p_n(\x)).
\end{align*}The alternative criteria proposed in \cite{Li12} are the following: \begin{align}\label{J_1}
\Ji(\x_{n+1}) & = \left.\mathbb{E}_n\left[\Big(\displaystyle\int_{\mathbb{X}}{\tau_{n+1}(\x)}^{\nicefrac{1}{2}}\textbf{P}_{\xx}(d\x)\Big)^2~\right|~\xx_{n+1}=\x_{n+1}\right],\\\label{J_2}
\Jii(\x_{n+1}) & = \left.\mathbb{E}_n\left[\Big(\displaystyle\int_{\mathbb{X}}{\nu_{n+1}(\x)}^{\nicefrac{1}{2}}\textbf{P}_{\xx}(d\x)\Big)^2~\right|~\xx_{n+1}=\x_{n+1}\right],\\\label{J_3}
\Jiii(\x_{n+1}) & = \left.\mathbb{E}_n\left[ \displaystyle\int_{\mathbb{X}}{\tau_{n+1}(\x)}\textbf{P}_{\xx}(d\x)~\right|~\xx_{n+1}=\x_{n+1}\right],\\\label{J_4} 
\Jiv(\x_{n+1}) & = \left.\mathbb{E}_n\left[\int_{\mathbb{X}}\nu_{n+1}(\x)\textbf{P}_{\xx}(d\x)~\right|~\xx_{n+1}=\x_{n+1}\right].
\end{align}
They are easier to estimate than $\JS$ because they are expressed  in terms of a single integral w.r.t. $\loiX$ (see \cite{Li12} for the implementation procedure with a Gauss-Hermite quadrature and a naive Monte Carlo method). 
Moreover, we precise that the Kriging mean $m_n$ is usually used as a deterministic interpolation model of function  $g$, so that  the quantity:\begin{align}\label{QSFDsDsDzd}\int_\E\mathbbm{1}_{m_n(\x)>T}\loiX(d\x),\end{align} is also an approximation of the failure probability $p$, which can easily be estimated by using a naive Monte Carlo method (see \cite{Dubourg13}, \cite{AKMCS} and \cite{Picheny10} for applications). As shown in \cite{Li12}, criteria $\Ji$ and $\Jiii$ actually refer to (\ref{QSFDsDsDzd}) for the estimation of $p$, while criteria $\Jii$ and $\Jiv$ refer to the expectation $\mu_n$, whose interest has been justified in Section \ref{First_Estimate}.

\begin{rem}\label{dldldld}
The strategies described above, which aim at enriching the design of experiments by adding points one by one, are usually called \enquote{1-step lookahead strategies}. In \cite{Chevalier14}, the authors show that it is possible to develop strategies, called  \enquote{$q$-step lookahead strategies}, which allow to add simultaneously $q\in\mathbb{N}^*$ points. Despite its interest, this approach is not addressed in this paper. 
\end{rem}

		\subsection{SUR strategy based on the convex order}\label{dudududu}\subsectionmark{}
 			\subsubsection{A new criterion}

Let us recall that the convex order inequality between the random variables $\Sn$ and $\Rn$ implies that:\begin{align*}\mbox{Var}[\Sn] \leq \mbox{Var}[\Rn], \quad\forall n \in\mathbb{N}^*.\end{align*}
Considering this inequality, the alternative strategy we suggest consists in determining the point $\x_{n+1}^*$ such that:\begin{align*}
\x_{n+1}^* = \argmin_{\x_{n+1}\in\E} \JR(\x_{n+1}),
\end{align*}where:\begin{align}\label{JR}
\JR(\x_{n+1}) = \mathbb{E}_n\big[\mbox{Var}[R_{n+1}]\mid \xx_{n+1} = \x_{n+1}\big].
\end{align}In Section \ref{Implementation}, we will explain how to estimate this criterion, the calculations being equivalent, in terms of algorithmic complexity, to those required for the estimation of the criteria $(J_{k, n})_{1\leq k \leq 4}$. Concerning the applications, we can directly consult the examples proposed in the Section \ref{Section7}. Before that, we propose in next sections to discuss the relevance of the criterion (\ref{JR}).

\begin{rem}
As mentioned in Section \ref{qhqye}, the quantity $H_{n+1}$ given in Equation (\ref{H_n}) refers to the posterior Bayes risk w.r.t. the quadratic loss function. It follows immediately, from Definition \ref{def_cx_order}, that one can derive an upper bound for any uncertainty measure $H_{n+1}$ based on a convex loss function. Thus, using the surrogate random variable $\Rn$ for $\Sn$, the scope of SUR strategies can be significantly extended.
\end{rem}

		\subsubsection{Comparison of criteria}

The fact that $\JR$ is always above $\JS$ is obvious and the following proposition shows in addition that it is locally closer than the functions $(\Jk)_{1\leq k\leq 4}$. 

\begin{pro}\label{Inequality_Criteria}
For all $k=1, \ldots, 4$, criteria $\JS$, $\JR$ and $\Jk$ satisfy:\emph{\begin{align}\label{<dsf<sdf}
\JS(\x)\leq \JR(\x)\leq \Jk(\x), \quad\forall\x\in\E. 
\end{align}}
\end{pro}
The partial order relation (\ref{<dsf<sdf}) simply means that, if one wants a local approximation of~$\JS$, then one must choose $\JR$. Nevertheless, we recognize that this does not guarantee that the approximation of the minimum of $\JS$ is better by considering $\JR$ rather than~$(\Jk)_{1\leq k\leq 4}$. To go further, we precise that the functions $\tau_n$ and $\nu_n$ can then be viewed as classification errors. Indeed, it is easy to verify that they reach their maximum when the membership probability  satisfies $p_n = \frac{1}{2}$. This corresponds to the dreaded  situation where the Kriging model is non-informative and fails to classify points in $\E$. In the opposite, they are equal to zero when $p_n = 0$ or $p_n=1$, which is a healthy situation.  
As a result, in order to provide information on the sub-domains of $\set$ where the Kriging model needs to be improved, the criteria $(\Jk)_{1\leq k\leq 4}$ tend to ideally select the next evaluation point among the set of points where classification errors are maximal. 
In the following, we show that criterion $\JR$ can be interpreted in the same way. For this purpose, we assume that the law of the random variable $p_n(\xx)$ is absolutely continuous, i.e. it admits a density. 

\begin{pro}\label{varR2_ClassError}The variance of $\Rn$ satisfies: \emph{\begin{align*}
\mbox{Var}[\Rn] = \displaystyle\int_\mathbb{X}\eta_n(\x)\textbf{P}_{\xx}(d\x),
\end{align*}}
where $\eta_{n}$ is a function of the membership probability  $p_n$, defined for all \emph{$\x\in\E$} by:
\emph{\begin{align*}
\eta_n(\x) = \big(1-p_n(\x)\big)\displaystyle\int_\E p_n(\y)\mathbbm{1}_{p_n(\y)\leq p_n(\x)}\textbf{P}_{\xx}(d\y)+p_n(\x)\displaystyle\int_{\E}\big(1-p_n(\y)\big)\mathbbm{1}_{p_n(\y)> p_n(\x)}\textbf{P}_{\xx}(d\y).
\end{align*}}
\end{pro}As for $\tau_n$ and $\nu_n$, 
the function $\eta_n$  is equal to zero when $p_n=0$ or $p_n = 1$. The following proposal shows that this function also has a maximum.
\begin{pro}\label{Max_etaFunc}Let us assume that the random variable \emph{$p_n(\xx)$} is absolutly continuous. Let $\eta_n:[0,1]\to\mathbb{R}$ be the function defined in Proposition \ref{varR2_ClassError} and $q_n^*\in [0,1]$ satisfying:\emph{\begin{align}\label{eq_etafunction}
{\int_\mathbb{X}}\mathbbm{1}_{p_n(\x) > q_n^*} \textbf{P}_{\xx}(d\x)= \mu_n \Leftrightarrow q_n^*= \mathbb{P}(\Rn>\mu_n).\end{align}}Then, $\eta_n$ has a global maximum on $[0,1]$ at $p_n = q_n^*$.
\end{pro}
In the end, we showed that the criterion $\JR$ can be rewritten:
\begin{align*}
\JR(\x_{n+1}) = \mathbb{E}_n\bigg[\int_{\E}\eta_{n+1}(\x)\loiX(d\x)\mid \xx_{n+1} = \x_{n+1}\bigg],
\end{align*} where the function $\eta_n$ can also be interpreted as a classification error, so that $\JR$ tends to ideally select the next evaluation location among the set of points with membership probability $p_n$ close to $q_n^*$. 

	\subsubsection{Interpretation based on random set theory}
The maximum value $q_n^*$ of $\eta_n$ depending on the law $\loiX$, it is not easy to analyse the point selection for the criterion $\JR$. Nevertheless, as shown in the sequel, the quantity $q_n^*$ can be linked to recent work on SUR strategies. More precisely, in order to understand the interest of Proposition \ref{Max_etaFunc} and derive an interpretation of criterion $\JR$, we introduce here some notions coming from the theory of random set (see \cite{TheoRandomSet} for an overview of random set theory, but also Chapter 2 of \cite{Chevalier13_Thesis} for applications in the setting of SUR strategies). In our context, the set $\Gamma_n $ defined by: \begin{align*}\Gamma_n= \lbrace\x\in\mathbb{X} : \xi_n(\x) > T\rbrace,\end{align*} is typically a random set. It satisfies $\loiX(\Gamma_n) = \Sn$ and consequently $\mathbb{E}[\loiX(\Gamma_n)] = \mu_n$. For all $\alpha\in(0,1)$, the $\alpha$-quantile of $\Gamma_n$ is the set $Q_{n}(\alpha)\subseteq\E$ satisfying: \begin{align*}
Q_{n}(\alpha) = \lbrace\x\in\mathbb{X} : p_n(\x) \geq \alpha\rbrace.\end{align*}  
In particular, the set $Q_n\left(\frac{1}{2}\right)$ is called  \enquote{Vorob'ev median} and $Q_n(q_n^*)$ \enquote{Vorob'ev expectation} (see \cite{TheoRandomSet}). According to Proposition \ref{Max_etaFunc}, we also have $\loiX(Q_n(q_n^*)) = \mu_n$. Let $A\Delta B = (A\cap B^c)\cup(A^c\cap B)$ denotes the symmetric difference between two sets $A$ and $B$. For any set $Q\subseteq\E$, we consider the distance $\mathbb{E}[\textbf{P}_{\xx}(\Gamma_n\Delta Q)]$ between $\Gamma_n$ and $Q$, which is called the \enquote{expected distance in measure}, e.g., in \cite{RandomSetsGauss}. The Vorob'ev median and expectation can be viewed as deterministic approximations of random set $\Gamma_n$ and be compared in terms of expected distance in measure. For the Vorob'ev expectation, it corresponds to the quantity:  \begin{align}\label{eterfetdre}
\mathbb{E}\big[\textbf{P}_{\xx}(\Gamma_n\Delta Q_n(q_n^*))\big] & = \mathbb{E}\left[\int_\set\mathbbm{1}_{\xi_n(\y)>T}\mathbbm{1}_{p_n(\y)\leq q_n^*}\loiX(d\y) + \int_\set\mathbbm{1}_{\xi_n(\y)\leq T}\mathbbm{1}_{p_n(\y)> q_n^*}\loiX(d\y)\right]  \nonumber \\ & = \int_\set p_n(\y)\mathbbm{1}_{p_n(\y)\leq q_n^*}\loiX(d\y) + \int_\set (1-p_n(\y))\mathbbm{1}_{p_n(\y)> q_n^*}\loiX(d\y), 
\end{align}which is called the \enquote{Vorob'ev deviation}. 
As mentioned in Chapter 2 of  \cite{TheoRandomSet}, the following proposition stands:
\begin{Pte}\label{jdybsy}
\begin{enumerate}
\item  For every $\alpha$-quantile $Q_n(\alpha)$ of the random set $\Gamma_n$, the Vorob'ev median $Q_n\left(\tfrac{1}{2}\right)$ satisfies:\emph{\begin{align*}\mathbb{E}\big[\loiX\left(\Gamma_n\Delta Q_n\left(\tfrac{1}{2}\right)\right)\big] \leq \mathbb{E}\big[\loiX \left(\Gamma_n\Delta Q_n(\alpha)\right)\big].\end{align*}}
\item For every set $Q$ satisfying \emph{$\loiX(Q) = \mu_n$}, the Vorob'ev expectation $Q_n(q_n^*)$ satisfies: \emph{\begin{align*}\mathbb{E}\big[\loiX \left(\Gamma_n\Delta Q_n(q_n^*)\right)\big] \leq \mathbb{E}\big[\loiX \left(\Gamma_n\Delta Q\right)\big].\end{align*}}
\end{enumerate}
\end{Pte}


Inequalities in Property \ref{jdybsy} implies that sets $Q_n\left(\tfrac{1}{2}\right)$ and $Q_n(q^*)$ satisfy:\begin{align*}\mathbb{E}\big[\loiX\big(\Gamma_n\Delta Q_n\left(\tfrac{1}{2}\right)\big)\big]\leq \mathbb{E}\big[\loiX \left(\Gamma_n\Delta Q_n(q_n^*)\right)\big],\end{align*}
meaning that, if one wants to approximate $\Gamma_n$ by a $\alpha$-quantile, then one must give priority to $Q_n\left(\tfrac{1}{2}\right)$. Otherwise, one can choose $Q_n(q^*)$, so that the volume of the approximation is equivalent to the average volume of $\Gamma_n$. 
By taking into account the remarks of the previous section, we conclude that criteria $(J_{n,k})_{1\leq k \leq 4}$ tend to preferentially select points in the neighbourhood of the boundary of $Q_n\left(\frac{1}{2}\right)$, which the best approximation of the random set $\Gamma_n$ in sense of the expected distance in measure, while the criterion $\JR$ preferentially selects points in the neighbourhood of the boundary of the set $Q_n(q^*)$, which is the best approximation of the random set $\Gamma_n$ among all sets with volume equal to the average volume of $\Gamma_n$. Note that, since we are in the first place interested in the failure probability $p$, estimated by measuring the average volume of $\Gamma_n$, it is not aberrant to use, as an approximation of $\Gamma_n$, an alternative set of equal volume on average. Moreover, we precise that authors  in \cite{Chevalier13_Thesis} also propose a criterion based on the Vorob'ev deviation given in Equation (\ref{eterfetdre}). Let us denote by $D_{Dev, n}$ this quantity and by $J_{Dev, n}$ the criterion based on it: 
\begin{align*}
J_{Dev,n}(\x_{n+1}) = \mathbb{E}_n\left[D_{Dev, n+1}\mid\xx_{n+1} = \x_{n+1}\right], \quad\forall\x_{n+1}\in\set.
\end{align*}Recent applications of this criterion can also be found in Chapter 5 of \cite{Azzimonti_Thesis}, as well as in \cite{elamri}. The authors use in particular to construct strategies that select several points simultaneously (see Remark \ref{dldldld}). Besides, by denoting $\eta_n^*$ the maximal value of function~$\eta_n$, i.e. when $p_n = q_n^*$, we can  verify that:\begin{align*}
\eta_n(\x)\leq \eta_n^* = \dfrac{D_{Dev, n}}{2}, \quad\forall\x\in\set,
\end{align*} meaning that the criterion $\JR$ satisfies:\begin{align*}
\JR(\x_{n+1})\leq \frac{J_{Dev,n}(\x_{n+1})}{2}, \quad\forall\x_{n+1}\in\set.%
\end{align*}Given this total order, it would be interesting to compare the performance of these two criteria on a concrete case. Indeed, it appears that $\JR$ is, once again, a better approximation of the reference criterion $\JS$. This paper do not study this approach in more detail, but it is a very interesting perspective considering recent work on this subject.

		\subsection{Implementation}\label{Implementation}

The implementation procedure for $\JR$ is basically the same as for $(J_{k,n})_{1\leq k\leq 4}$, as explained in Section 3.3 and Section 3.4 of \cite{Li12}. Typically, since we have in the case of a Gaussian Kriging model: \begin{align*}\JR(\x_{n+1}) & = \mathbb{E}_n\big[\mbox{Var}_{}[R_{n+1}]\mid\xx_{n+1} = \x_{n+1}\big] \\
& = \mathbb{E}_n\big[\mbox{Var}_{}[R_{n+1}(\x_{n+1}, \xi(\x_{n+1}))]\big] \\
& = \int_{\mathbb{R}}\mbox{Var}_{}[R_{n+1}(\x_{n+1}, z)]\dfrac{1}{\sigma_n(\x_{n+1})\sqrt{2\pi}}\mathrm{e}^{-\frac{1}{2}\left(\frac{z - m_n(\x_{n+1})}{\sigma_n(\x_{n+1})}\right)^2}dz,
\end{align*}
we use a Gauss-Hermite quadrature to approximate this integral:
\begin{align*}
\JR(\x_{n+1}) & \approx \frac{1}{\sqrt{\pi}}\sum_{q=1}^Q w_q\mbox{Var}[R_{n+1}(\x_{n+1}, m_n(\x_{n+1})+\sqrt{2}u_q\sigma_n(\x_{n+1}))],
\end{align*}
where $(w_q)_{1\leq q\leq Q}$ and $(u_q)_{1\leq q\leq Q}$ stand for the quadrature weights and the quadrature points.\medskip

Moreover, to approximate the variance of $\Rn$, we do the following. Let~$(\xx_i)_{1\leq i \leq N}$ be an i.i.d. sample with distribution $\loiX$ and assume that, among probabilities $(p_n(\xx_i))_{1\leq i \leq N}$, there are $N'$ distinct. Their increasing reordering is the following: $$0\leq p_n^{(1)}\leq \ldots \leq p_n^{(N')}\leq 1,$$ where $p_n^{(1)}= \underset{1\leq i \leq N}\min(p_n(\xx_i))$ and $p_n^{(N')} =\underset{1\leq i \leq N}\max(p_n(\xx_i))$. For all~$1\leq i\leq N'$, we define~$l_i$\\

 as the number of occurrences of the probability $p_n^{(i)}$, so that $\sum_{i=1}^{N'}l_i = N$. In addition, we introduce the notation $n_i = N - \sum_{j=1}^{i}l_j$, where $n_{N'} = 0$ and $n_{N'-1}=l_{N'}$. 
As a result, we get:\begin{align*}
\mbox{Var}[\Rn] & = \int_{\E^2}\min\big(p_n(\x_1), p_n(\x_2)\big)\loiX(d\x_1)\loiX(d\x_2) - \left(\int_\E p_n(\x)\loiX(d\x)\right)^2\\
& \approx \dfrac{1}{N^2}\sum_{i=1}^{N'}l_ip_n^{(i)}(l_i+2n_i) - \left(\dfrac{1}{N}\sum_{i=1}^{N'}l_ip_n^{(i)}\right)^2.
\end{align*}
Note that, for $N'=N$, the latter is equal to the Monte Carlo estimator (\ref{lodkdh}).

\section{Numerical experiments}\label{Section7}
	\subsection{A one-dimensional example}
		
\begin{figure}[h!]
\begin{center}\hspace{0.1cm} 
\begin{minipage}[b]{1\linewidth}
\includegraphics[width=7cm, height = 4.5cm]{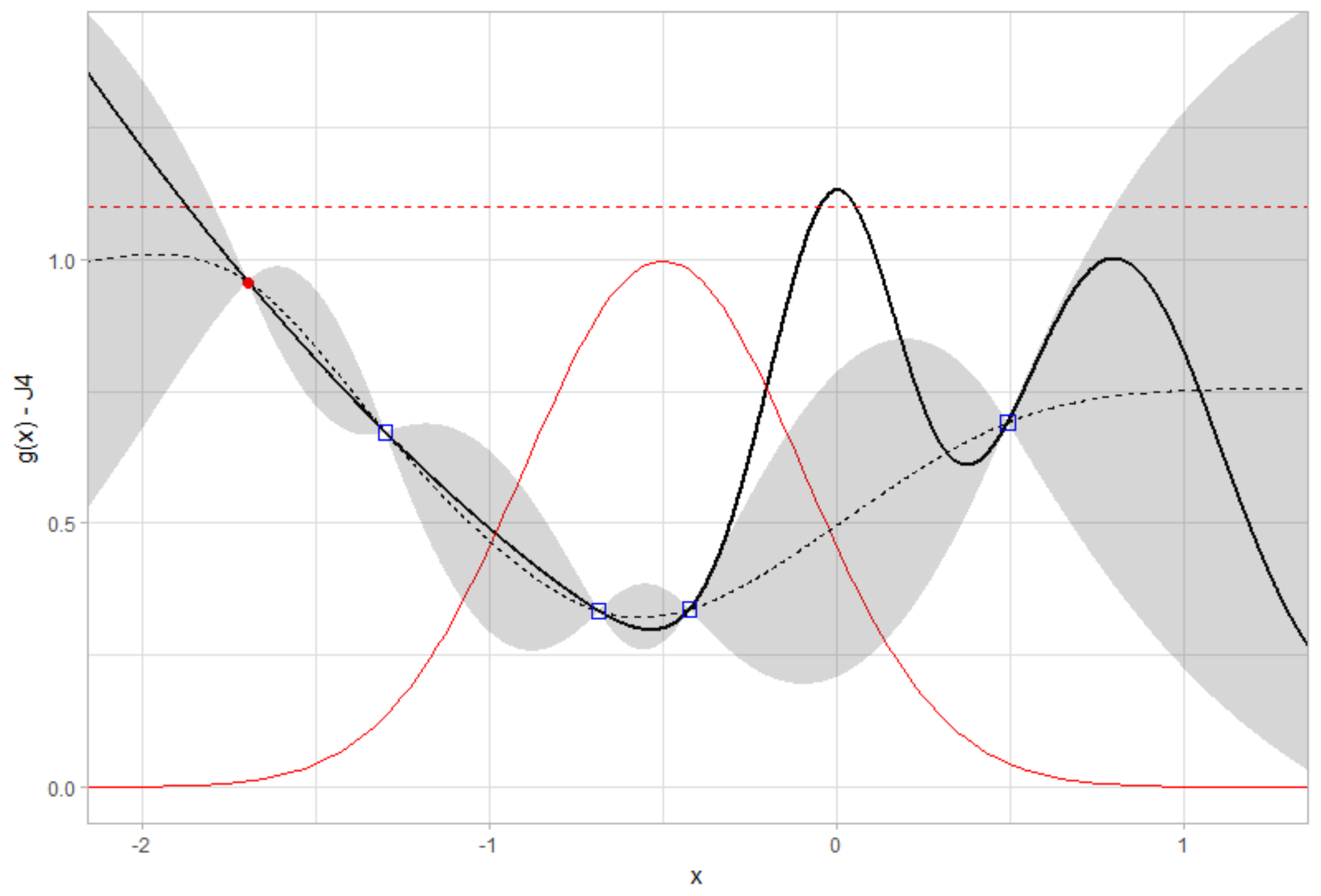}
\includegraphics[width=7cm, height = 4.5cm]{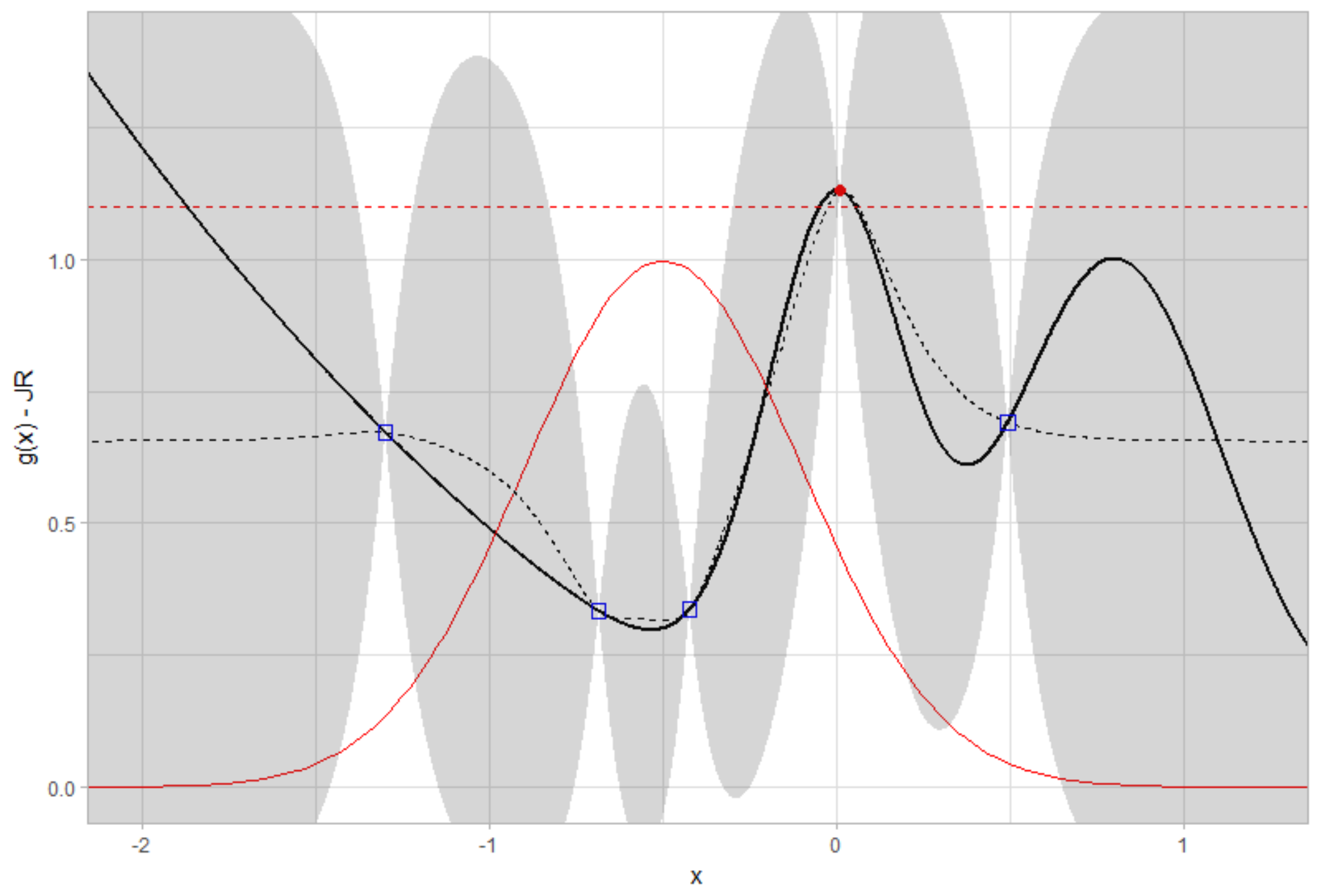}
\end{minipage}
\begin{minipage}[b]{1\linewidth}
\includegraphics[width=7cm, height = 4.5cm]{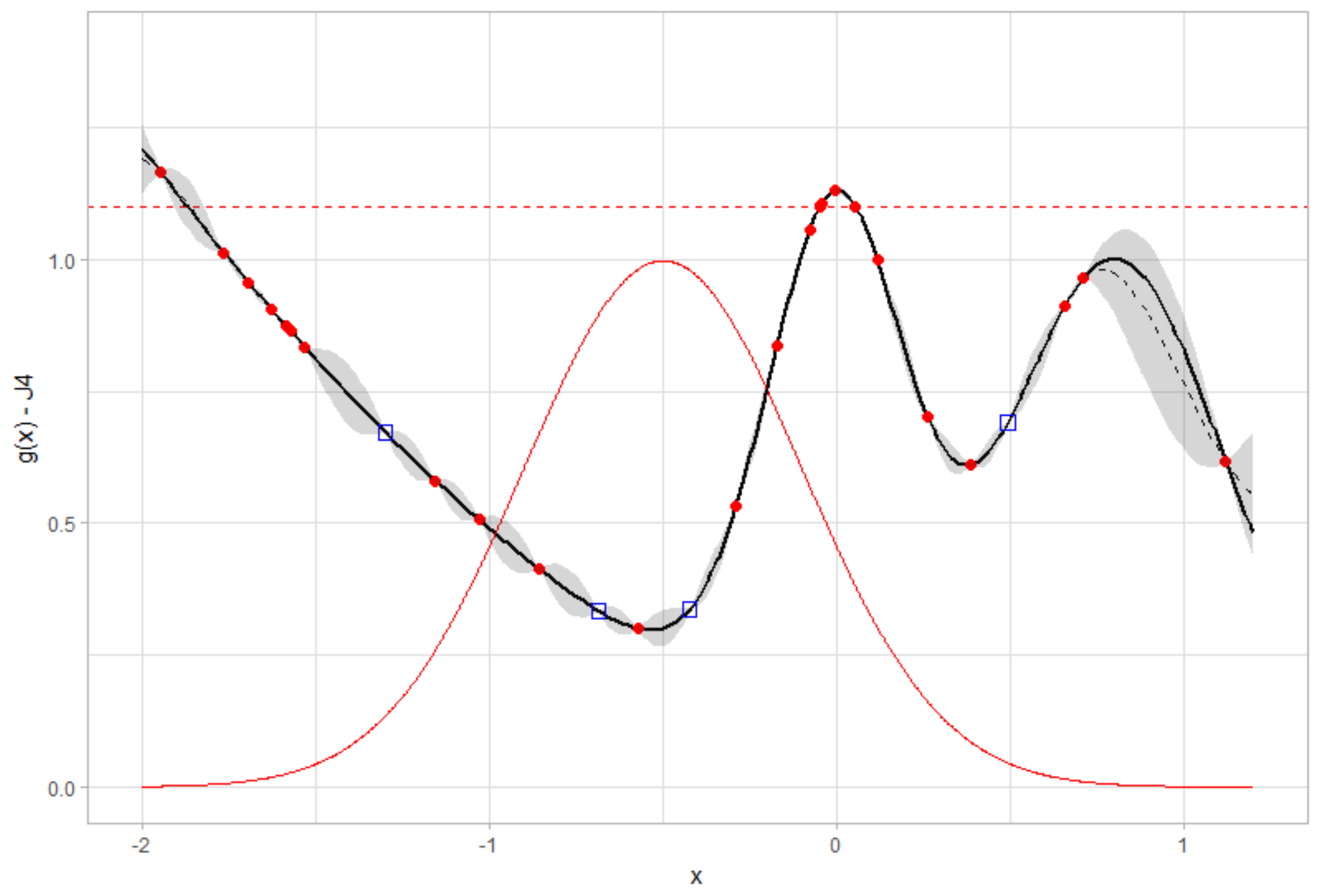}
\includegraphics[width=7cm, height = 4.5cm]{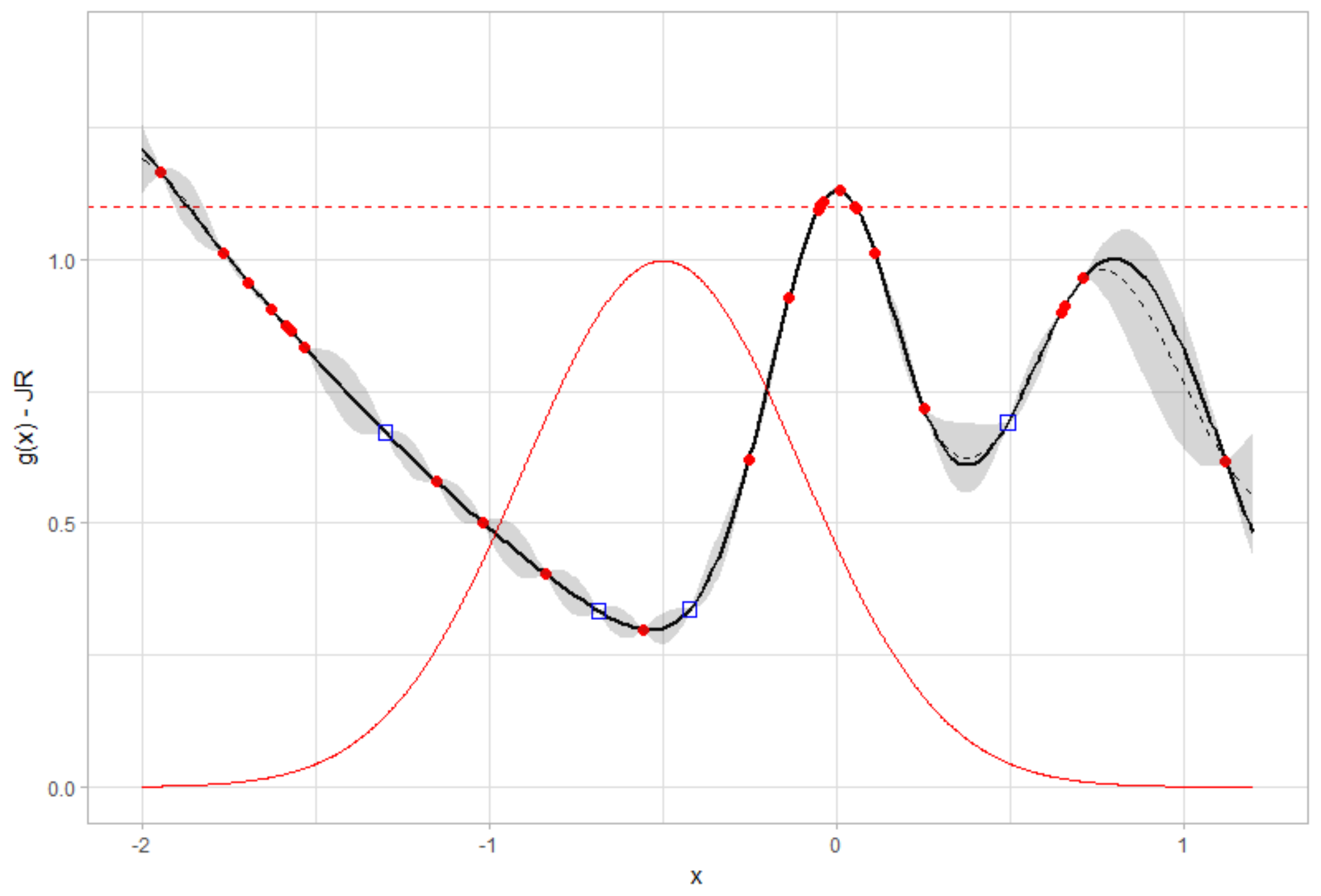}
\end{minipage}\vspace{-0.5cm}
\caption{{\small Illustration of SUR strategies based on criteria $J_{4,n}$ (left) and $\JR$ (right). Top: First iteration. Bottom: Last iteration. Function $g$ (black plain line); threshold $T$ (red dashed line); initial experiments (squares); the new ones (circles); mean function $m_n$ of the Kriging model (dashed curve); $95\%$ confidence intervals given in Equation (\ref{kdolemmmmm}) (shaded area); density of $\loiX$ (red plain line).}}
\label{g_Ex1}
\end{center}
\end{figure}

Here, we consider the efficiency of our estimation procedure on a simple artificial example, that is a modified version of the one proposed in \cite{Li12}. 
It consists in estimating the probability $p = \mathbb{P}\left(g(\xx)>T\right)$, where $T=1.1$, $\xx$ is a random variable with distribution $\textbf{P}_{\xx} = \mathcal{N}(-0.5, 0.4^2)$ and $g: \mathbb{X} = \mathbb{R} \rightarrow  \mathbb{R}^{+}$ is such that:
\begin{align*}g(\x)=(0.4\x-0.3)^2+e^{-11.534|\x|^{1.95}}+e^{-5(\x-0.8)^2}.\end{align*}We known in advance that $p = 4.643\cdot 10^{-2}$. 
We propose to evaluate our approach with~3 different design sampling methods: an LHS-maximin method that generates a design of 30 experiments (see, e.g., \cite{Dette}, \cite{plans_LHS} and \cite{lhs_stein}) and two sequential method respectively based on the use of criteria $\Jiv$ and $\JS$. They take as an input a design of size $n=4$ (see Figure \ref{g_Ex1}).

\begin{figure}[h!]
\begin{center}\hspace{0.1cm} 
\begin{minipage}[b]{1\linewidth}
\includegraphics[width=7cm, height = 4.5cm]{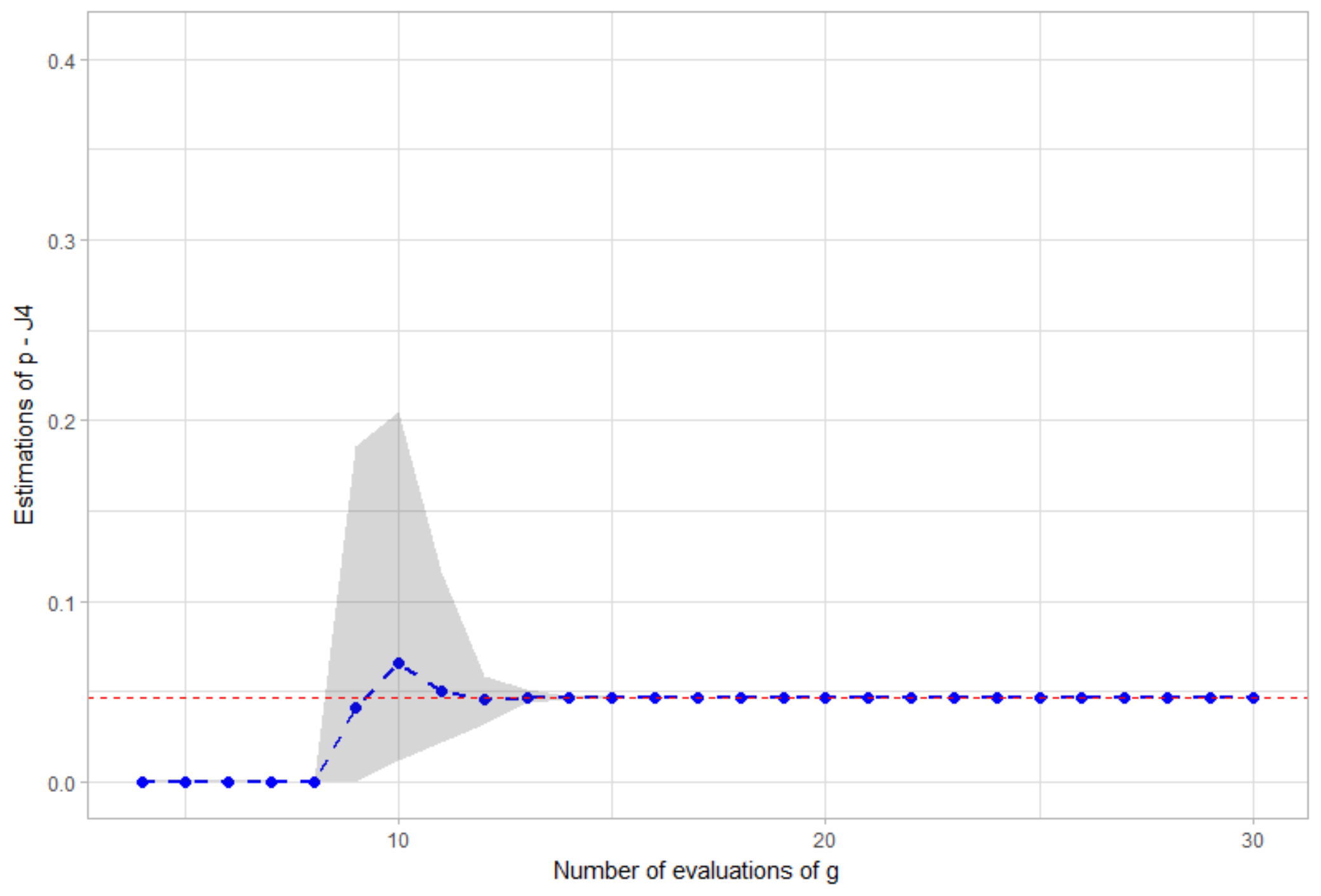}
\includegraphics[width=7cm, height = 4.5cm]{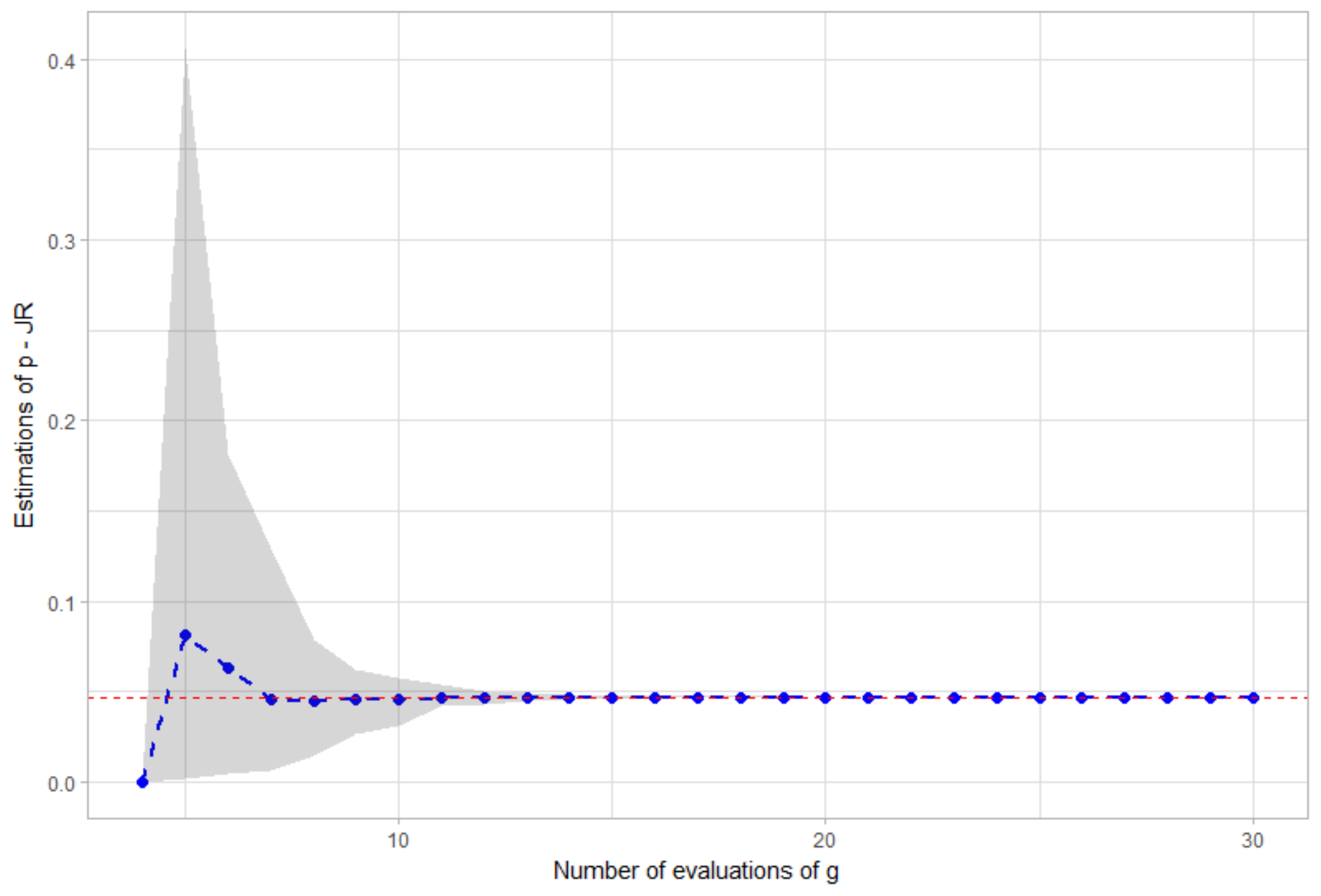}
\end{minipage}\vspace{-0.5cm}
\caption{{\small Illustration of performances of SUR strategies based on criteria $\Jiv$ (left) and $\JR$ (right). True failure probability $p~=~4.643\cdot~10^{-2}$ (red dashed line) ; successive estimates of $p$ corresponding to the estimator (\ref{lodkdh}) (circles); successive estimates of a credible interval at level $95\%$ given in Equation (\ref{Intervalle_cred_cx}), with $\beta = \frac{1}{2}$ (shaded area).}}
\label{est_Ex1}
\end{center}
\end{figure}

\begin{figure}[h]
\begin{center}\hspace{0.1cm} 
\begin{minipage}[b]{1\linewidth} 
\includegraphics[width=7cm, height = 4.5cm]{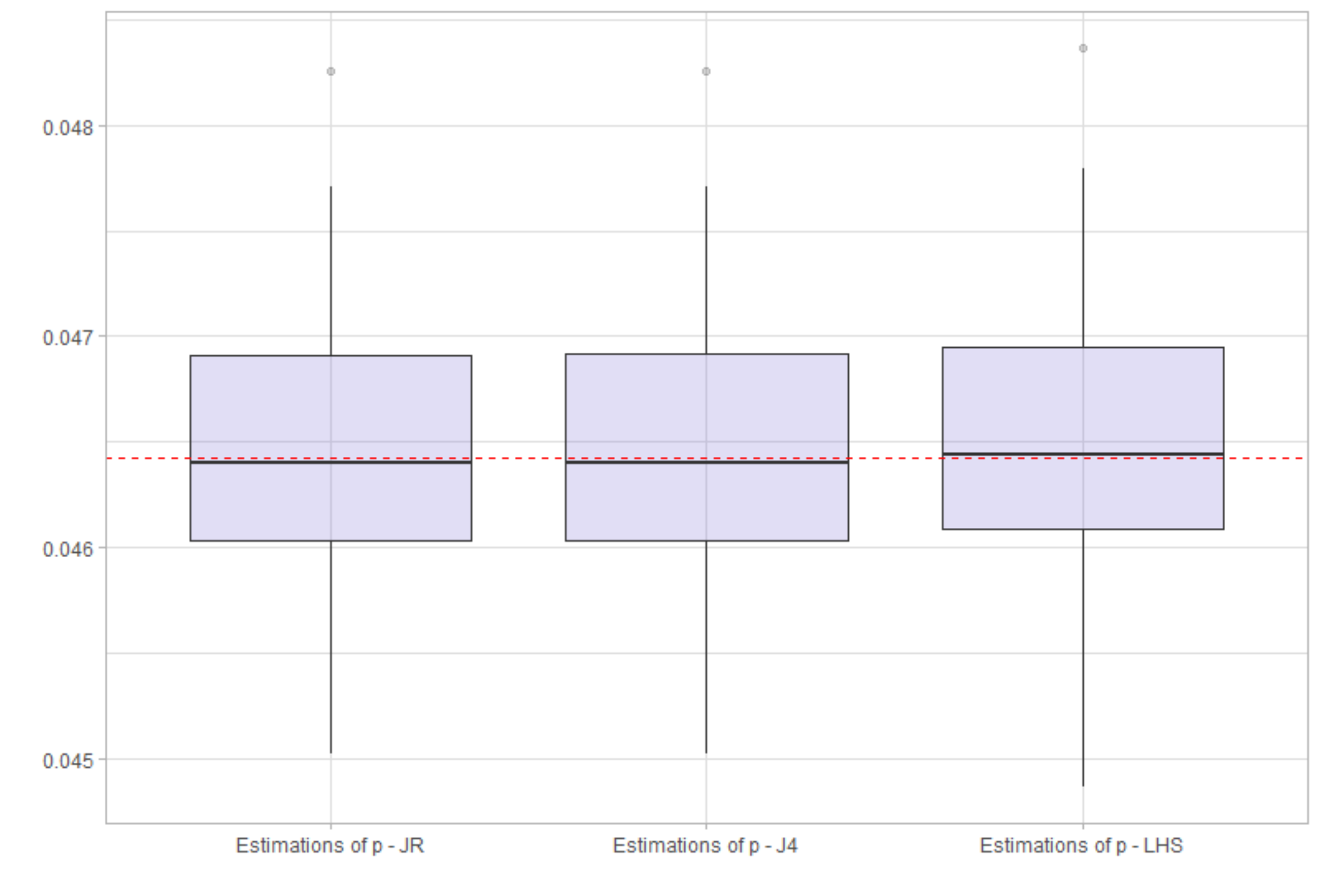}
\includegraphics[width=7cm, height = 4.5cm]{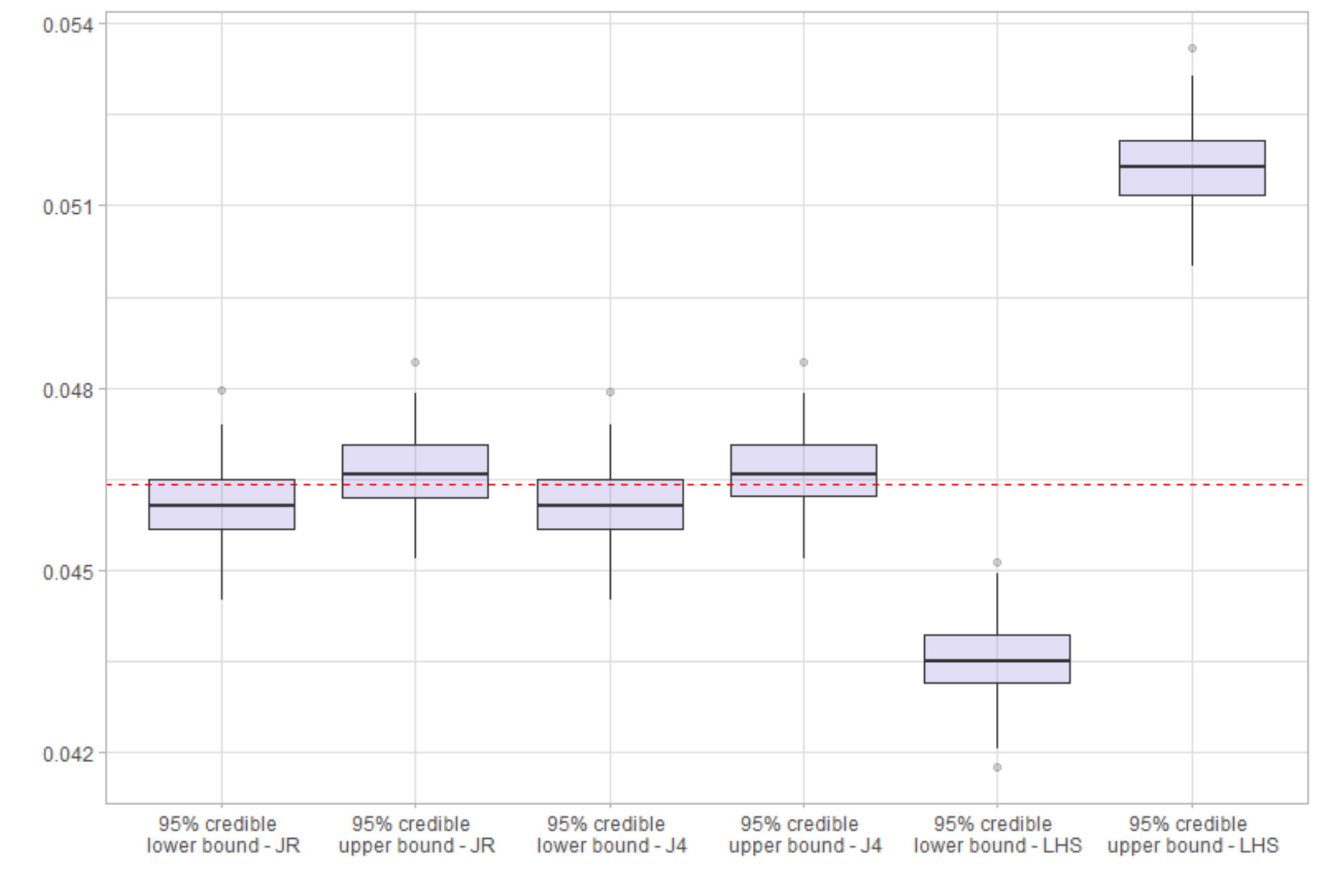}
\end{minipage}\vspace{-0.5cm}
\caption{{\small Sensitivity study to sample size of the Monte Carlo method. Left \& right: True probability $p = 4.643\cdot 10^{-2}$ (red dashed line). Left: Boxplots for the estimation of $p$ obtained with 100 Monte Carlo simulations. Right: Boxplots for the estimation of a credible interval at level $95\%$ obtained with the same 100 Monte Carlo simulations.}}\label{Boxplot}
\end{center}
\end{figure}

 We select $\Jiv$ because it seems to be the criterion which is the most used in practice (see, e.g., \cite{KrigInv}). 
In order to build the Kriging model, we use the function \enquote{km} of the {R} package called \enquote{DiceKriging} (see \cite{DiceKriging}). By using a Leave-One-Out cross validation method, we choose a constant mean and a Mat\'ern  covariance function. For the sequential strategies, we update the Kriging model at each step by taking account of the new point and its evaluation. In Figure \ref{g_Ex1}, the initial design is represented (blue squares) and the successive choices of points according to the criterion considered can be compared (red points). Starting from a number of simulations, the two criteria select essentially the same points and the evaluations focus on the neighbourhood of the threshold $T$ (red dashed line). In Figure \ref{est_Ex1}, we observe the progress of each SUR strategy in terms of estimations. The estimates of $p$ (blue points) corresponding to (\ref{lodkdh}) and the related credible intervals (gray area), which are derived from (\ref{Intervalle_cred_cx}) with $\beta = \frac{1}{2}$, are all obtained performing a naive Monte Carlo method. 
They must be compared with the true value of $p$ (red horizontal dashed line). 
 According to Figure \ref{est_Ex1}, the estimations provided by $\JR$ and $\Jiv$ are in both cases consistent with the true value of $p$. Moreover, as can be seen that the credibility intervals tends to decrease as the number of observations increases. When $n=30$ observations of the function $g$ are available for each design sampling methods, we use an i.d.d. sample w.r.t.~$\textbf{P}_{\xx}$ of size $N = 10^4$ to provide (with the \enquote{predict} function from \enquote{DiceKriging} package) an estimation of $p$ and  a $95\%$-credible interval. We run this operation 100 times in order to get boxplots as in Figure \ref{Boxplot}. In Table \ref{Values_Ex1}, the corresponding average values are reported. We see that both sequential strategies provide very satisfactory similar results. The LHS design also leads to a good estimation of $p$, but the corresponding credible interval is significantly larger, meaning that the estimation of $p$ is less reliable. Besides, we precise that the Markov's bounds of Proposition \ref{Prop_CredibleInt_Markov} are, in this example, not informative since they are always equal to 0 and 1. 

\begin{table}[h!]
\begin{center}
\renewcommand{\arraystretch}
{1.2}
\begin{tabular}{|c|c|c|c|}
\hline
 & {\small Estimation of $p$} & {\small $95\%$ credible lower bound} &  {\small $95\%$ credible upper bound} \\
\hline
{\small $\Jiv$} & {\small $4.644\cdot10^{-2}$}   &  {\small $4.608\cdot10^{-2}$} & {\small $4.663\cdot10^{-2}$} \\
\hline
{\small $\JR$} & {\small $4.645\cdot10^{-2}$}  & {\small $4.608\cdot10^{-2}$}  & {\small $4.663\cdot10^{-2}$}\\
\hline
{\small LHS} &  {\small $4.649\cdot10^{-2}$} &  {\small $ 4.352\cdot10^{-2}$}  & {\small $5.163\cdot10^{-2}$ } \\
\hline
\end{tabular}
\caption{{\small Estimations and credible bounds at level $95\%$ for the target probability $p=4.642\cdot10^{-2}$, defined as the average boxplot values of Figure \ref{Boxplot}.}}
\label{Values_Ex1}
\end{center}
\end{table}

\subsection{An industrial case study}

		\subsubsection{Description of the real case}

The estimation methods described in this paper have been tested on a real case of the company STMicroelectronics. It
concerns the study of an electronic component called duplexer. It is a device used to filter a signal over several frequency bands in order to extract and isolate different part of it. As STMicroelectronics covers the mobile telephony market, the duplexer studied operates in the radiofrequencies ($f$ $\thicksim$ 1 GHz). The main signal transmission channel is divided into two distinct channels. Each channel has a bandpass filter to select certain frequency components of the signal. The semi-conductor manufacturing process used to make bandpass filter is a succession of complex operations, difficult to maintain constant over time: deposition of metal layers by electrolysis, plasma etching, photolithography\ldots The resulting technology, called RLC06A, can then be seen as a stack of insulating and conductive layers. The parameters subject to variations are typically the thicknesses of these deposits. For the duplexer studied, exactly 4 deposit thicknesses are influential, in the sense that their variations significantly impact filter performance. It is possible to model the natural variability of the industrial process by associating to these quantities a probability distribution, so that they become random variables. They are noted $X_1, \ldots,X_4$ and their distributions are extracted by from online measurements (i.e. directly on the production lines) of test patterns. These can be properly approximated by normal distributions. The parameters of each of them are given in~Table~\ref{Tableau_Cas_Indus1}.

\begin{table}[h!]
\begin{center}
\renewcommand{\arraystretch}
{1.2}
\begin{tabular}{|c|c|c|c|c|c|c|}
  \hline
{\small Factor}  & {\small Name}   & {\small Unit}  & {\small Minimal value} & {\small Maximal value}  & {\small Distribution}  & {\small Position and scale parameters} \\
  \hline
{\small $X_1$} &  {\small BCB1}   & {\small $\upmu$m} &{\small 2.7504}  & {\small 5.1011}  & {\small Gaussian}  & {\small 3.3477~;~\small 0.19108} \\
  \hline
{\small $X_2$} &  {\small Capa2}  & {\small $\upmu$m}   & {\small 160.41}  & {\small 188.31}  & {\small Gaussian} & {\small 174.31~;~\small 1.6831} \\
  \hline
{\small $X_3$} &  {\small Meta1b} & {\small $\upmu$m} & {\small 0.62241} & {\small 0.90364} & {\small Gaussian} & {\small 0.73389~;~\small 0.03193} \\
  \hline
{\small $X_4$} &  {\small Meta2}  & {\small $\upmu$m} & {\small 4.9964}  & {\small 7.4946}  & {\small Gaussian} & {\small 6.1457~;~\small 0.2678} \\
  \hline
\end{tabular}
\caption{\small{Distribution of the 4 thicknesses of deposits that impact the duplexer response.}}
\label{Tableau_Cas_Indus1}
\end{center}
\end{table}
In the context of high-frequency electronics, it is shown that the duplexer can be completely characterized by its dispersion matrix. This is a $3\times $3 matrix, which gives for each input-output of the device the proportion of the signal transmitted or reflected. Each parameter depends on the frequency. The customer's  specifications are defined on these frequency characteristics: it gives a template that the frequency response of the duplexer must respect. This template does not cover the entire frequency range but only certain bands. Here, this translates into 12 responses (vectorial ouput), each of which is an unknown function of the 4 input variables given in the table \ref{Tableau_Cas_Indus1}. Each response is associated with a frequency range and has its own constraint, expressed in decibels.   All the responses, their characteristics and associated constraints  can be found in details in the thesis \cite{lucie_Thesis}, where this industrial case is also treated. 
Numerical simulations are performed using commercial software {HFSS} (High Frequency Structural Simulator, {ANSYS}$^{\mbox{\scriptsize{\textregistered}}}$ Electronics Desktop, Release 17.2, {ANSYS} Inc), which provides results for any virtual configuration of the product. 
 
		\subsubsection{Formalism et preliminary study}
Here, a numerical simulation (finite element analysis) via the {HFSS} software requires the execution of a calculation code, represented by a function $g : \x\in\set \mapsto g(\x) = (g_1(\x),\ldots,g_{12}(\x))\in\mathbb{R}^{12}$, where $\set\subseteq\mathbb{R}^4$ is defined is Table \ref{Tableau_Cas_Indus1}. 
The product is considered defective if at least one of the response is defective, so that the probability of failure $p$ is written:\begin{align}\label{ppa}
 p = \mathbb{P}\left(\bigcup_{j=1}^{12}g_j(\xx) \geq T_j\right),
\end{align} where $\xx = (X_1, \ldots, X_4)$ is given in Table \ref{Tableau_Cas_Indus1}. From a formal point of view, the methods developed in this paper can still be applied, with the construction of a multi-dimensional random process. Indeed, there are Kriging models adapted to the case where the output is vectorial, called co-Kriging models (see \cite{co-krig}, \cite{co-krig2} and the thesis \cite{LeGratiet_Thesis}). Despite their interest, they are not integrated here because they introduces considerable complexity for practical implementation. 
We explain in the next section how to deal simply with the multi-dimensional case. Before that, we precise that a preliminary study of this industrial case was carried out performing $N=1000$ simulations using the software {HFSS}, the idea being to estimate the probability of failure using a naive Monte Carlo method. Given a sample $(\xx_i)_{1\leq i \leq N}$ i.i.d. distributed according to the laws of the input factors (see Table \ref{Tableau_Cas_Indus1}), 
the estimator $\widehat{p}_N$ of 
probability (\ref{ppa}) is:\begin{align*}\widehat{p}_N 
= \frac{1}{N}\sum_{i=1}^{N}\mathbbm{1}_{\bigcup_{j=1}^{12}g_j(\xx_i)\geq T_j}.
\end{align*}Thus, the probability of failure is approximately $6.7\times10^{-2}$ and, by the central limit theorem, a  $95\%$-confidence interval is $[4.5\times10^{-2},~ 8.8\times10^{-2}]$. These simulation results will be used as a reference to validate the results obtained. 

\subsubsection{Strategy for vector output}		

Suppose that the function $g$ is already known at the points of a design of experiments $(\x_i)_{1\leq i \leq n}$. For each responses, these observations are used to construct a Gaussian process model. For every $j=1, \ldots, 12$, we note $\xi_{n, j}$ the Gaussian process conditioned on observations corresponding to the $j$-th response.  For all $\x\in\set$, we define:\begin{align}\label{kidheutj}
p_n(\x) = \mathbb{P}\left(\displaystyle\bigcup_{j=1}^{12}\xi_{n, j}(\x)\geq T_j\right).
\end{align}
 For all  $j=1, \ldots, 12$, by denoting $p_{n,j}(\x) = \mathbb{P}\left(\xi_{n, j}(\x)\geq T_j\right)$, we also define: \begin{align*}
p_n^+(\x) = \min\left(1, \sum_{j=1}^{12}p_{n,j}(\x)\right).
\end{align*}Since the quantity (\ref{kidheutj}) is here inaccessible, we bound it from above. Indeed, the union boundary implies that:\begin{align}\label{efqqefqfaeafaeff}
p_n(\x)\leq  p_n^+(\x), \quad\forall\x\in\set.
\end{align}If all models are sufficiently informed, the probabilities $(p_{n, j}(\x))_{1\leq j \leq 12}$ are close to 0 or 1, and this does not induce any significant error. As a result, the estimation of probability~(\ref{ppa}), which is considered in practice, has the expression $
\int_\set p_n^+(\x)\loiX(d\x)
$. Note that doing so, we majorize the quantity of interest (\ref{Mean_R1}), unknown in the multivariate case. Such a precaution is key in risk assessment. Moreover, the $95\%$-credibility interval considered is derived from Proposition \ref{dsjsdjzDJsef}, with $\beta = \frac{1}{2}$, and defined by:
\begin{align}\label{mmmmmdyyybs}
\left[1-\int_\set\min\left(1, \frac{1-p_n^+(\x)}{0.025}\right)\loiX(d\x)\quad,\quad\int_\set\min\left(1, \frac{p_n^+(\x)}{0.025}\right)\loiX(d\x)\right].
\end{align} All the integrals above are evaluated numerically  applying a naive Monte-Carlo method. Concerning the sequential {SUR} strategy based on the variance of $\Rn$, we proceed as follows. For each response, a Gaussian process model is constructed and we derive the variance of the random variable $\Rn$. To determine the point $\x_{n+1}$ where to perform the next $g$ evaluation, the sampling criterion $\JR$, which corresponds to the response for which the variance of $\Rn$ is the highest, is optimized. This procedure simply consists in selecting the response for which the uncertainty about the result is the highest. Then, the selected point is added to the existing design of experiments and, finally, the values of the 12 responses are calculated. 

		\subsubsection{Results}
	
\begin{figure}
\begin{minipage}[b]{1\linewidth}
\begin{center}
\includegraphics[width=8cm, height = 5.5cm]{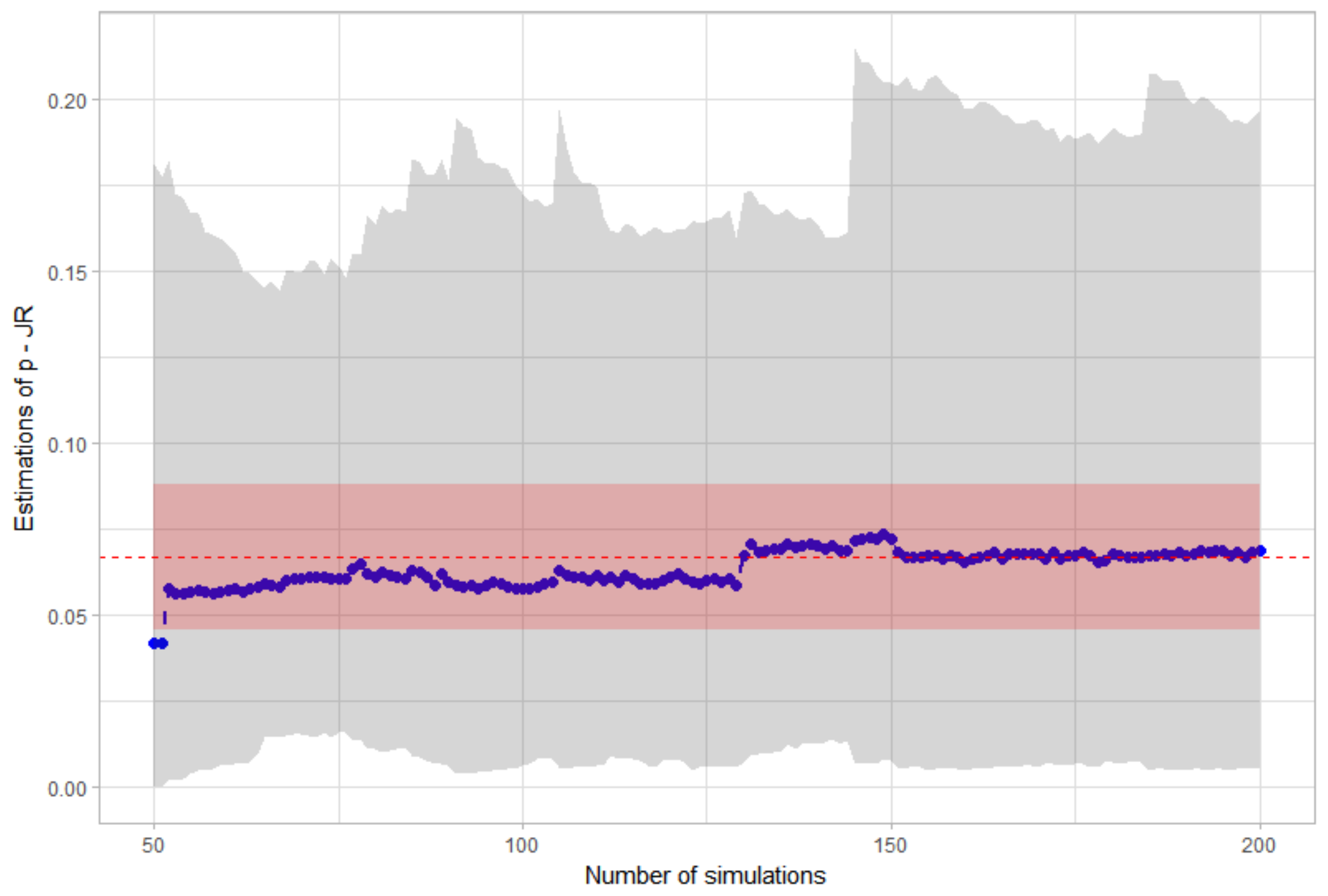}
\end{center}
\end{minipage}\vspace{-0.5cm}
\caption{{\small Successive estimations of of the failure probability $p$ (horizontal red dashed line); $95\%$-confidence intervals with $N=1000$ simulations (red area); estimations of $p$ (blue points) and $95\%$-credible intervals (gray area) with $n=50, \ldots, 200$ simulations.}}
\label{ggg}
\end{figure}

For each response, the Kriging model has a constant trend and an isotropic Mat\'ern correlation function with parameter $\frac{3}{2}$. These choices have been validated by a cross-validation method (for more information, see \cite{lucie_Thesis}). Our initial design of experiments is a {LHS} design of size~$n=50$. We add 150 points by applying the {SUR} strategy based on the $\JR$ criterion. However, instead of optimizing the criterion on the whole set~$\set$, and in order to avoid additional simulations, we determine the point to add to the design of experiments from the set of points $(\xx_i)_{1\leq i \leq N}$ whose values of $g$ are already available. The results we obtained are given in Table \ref{Tableau_kdkdjdd}, and the successive estimates (blue points) are represented in Figure \ref{ggg}. The estimated value of $p$ by the Monte Carlo method is represented by the red dashed line: this is our reference measure, although it is not the true value of $p$. With only $n=150$ simulations, the estimate of the failure probability obtained by applying the {SUR} strategy based on $\JR$ is equivalent to this latter. From this point of view, the results are satisfactory. Nevertheless, we recognize that the credibility intervals obtained (gray area) are not very informative in this case, except to highlight a modeling issue: the addition of data here does not seem to result in a significant selection of process trajectories compatible with them. It is possible that selecting the points to add to the design of experiments from the set $(\xx_i)_{1\leq i \leq N}$ decreases the effectiveness of the strategy, or that the correlation function (and therefore the interpolation basis) chosen is not adapted in the end for at least one response. 

\begin{table}
\begin{center}\renewcommand{\arraystretch}
{1.5}
\begin{tabular}{|c|c|c|c|}
  \hline
   & \small{Mean value of $\Rn$} & \small{Variance of $\Rn$} & \small{$95\%$-credible interval}\\
  \hline
\small{$J_n^{\Rn}$ ; $n = 200$} & 6.9$\times10^{-2}$
& \small{1$\times10^{-3}$} 
& \small{$\big[$1.3$\times10^{-3}$;~1.9$\times10^{-2}\big]$ }
\\
  \hline
\end{tabular}
\caption{\small{Final estimation of the failure probability $p$ via a {SUR} strategy based on the $\JR$ criterion.}}
\label{Tableau_kdkdjdd}
\end{center}
\end{table}

\begin{figure}[h]
\begin{center}\hspace{0.1cm} 
\begin{minipage}[b]{1\linewidth} 
\includegraphics[width=7cm, height = 4.5cm]{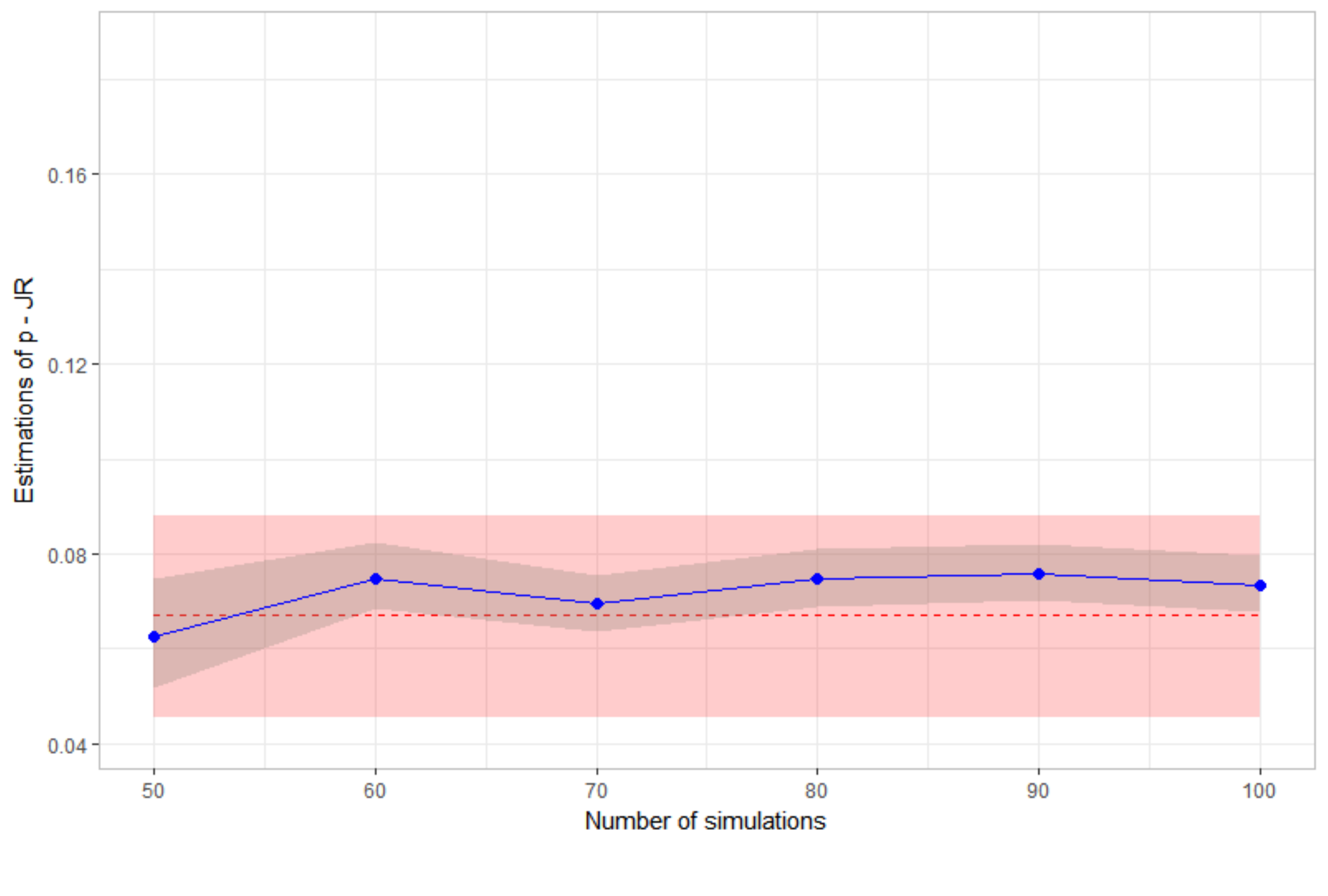}
\includegraphics[width=7cm, height = 4.5cm]{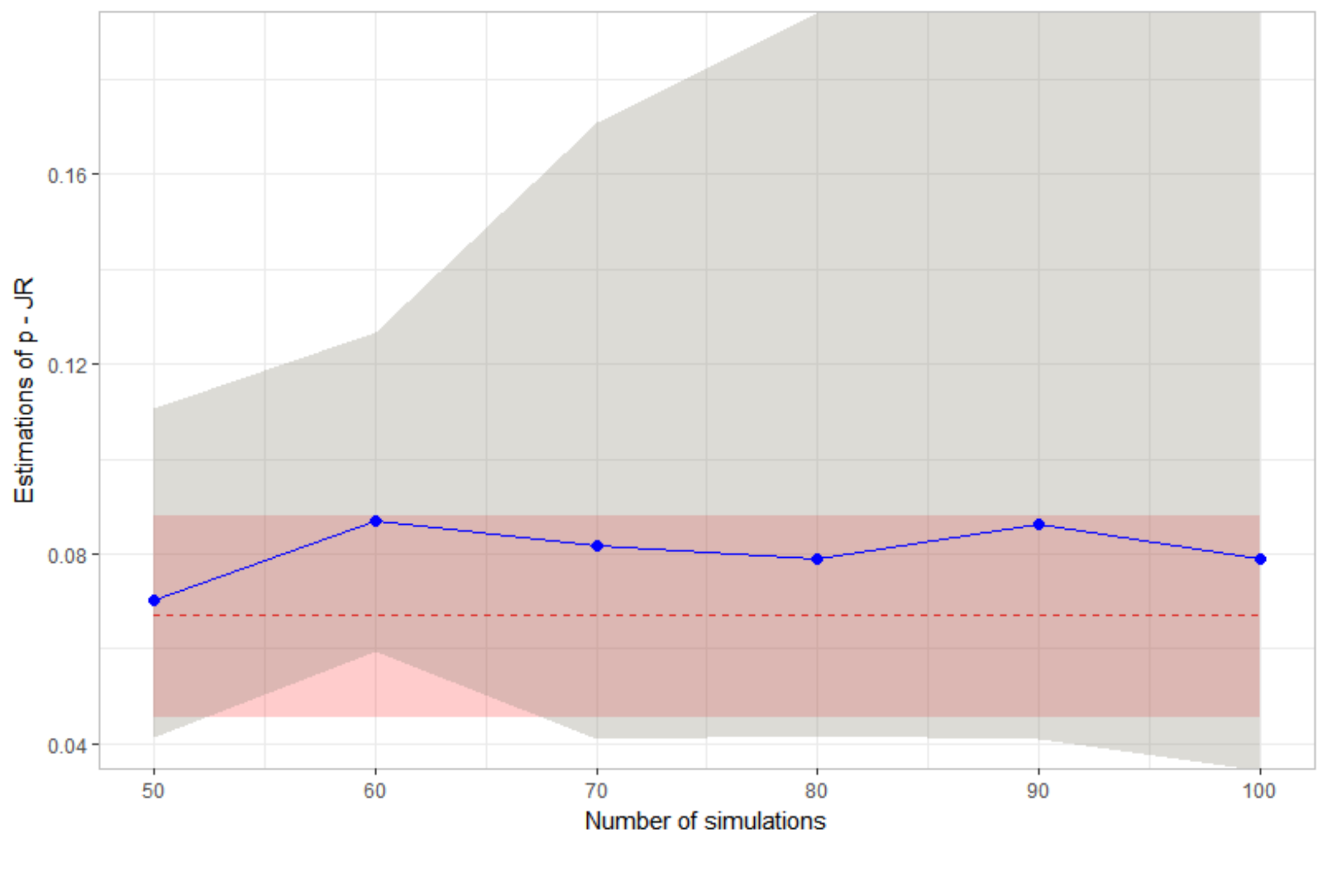}
\end{minipage}\vspace{-0.5cm}
\caption{{\small Successive estimations, with consideration of the derivative (left) and without consideration of the derivative (right), of the failure probability $p$ (horizontal red dashed line); $95~\%$-confidence intervals with $N=1000$ simulations (red area); estimations of $p$ (blue points) and $95\%$-credible intervals (gray area) with $n=50, \ldots, 200$ simulations.\vspace{-0.5cm}}}\label{ggg1}
\end{center}
\end{figure}

		\subsubsection{Study with consideration of the derivative}

When performing the numerical simulations, it is possible in our case to know the values of the derivatives at the measured points. Indeed, the {HFSS} simulator can also provide the different partial derivatives of $g$ according to each dimension, i.e. the gradient. For each enhanced output, we can therefore construct a Kriging model that takes into account the gradient at the evaluation points. For more information on taking the value of the derivative into account in Gaussian process modeling, see, e.g. \cite{deriv_krig} and \cite{krig_deriv2}. Note that, when including derivatives, it becomes fundamental to implement a robust method to set the values of the hyper-parameters as well as to invert the covariance matrix. Here, we chose again a constant trend and an isotropic Mat\'ern correlation function (see \cite{lucie_Thesis} for a discussion on this model choice). The results obtained by incorporating the observed gradients are given in Figure \ref{ggg1} (left). They should be compared with those given in Figure \ref{ggg1} (right), where the Kriging model is constructed without taking into account the values of the derivatives. We treated these both cases with the same initial design of type {LHS} consisting of $n=50$ points. Note that, this time, the candidate points are not chosen from the sample $(\xx_i)_{1\leq i\leq N}$ to build the experimental design, since the simulations are carried out at points  identified by the criterion $\JR$. 
In Figure \ref{ggg1}, we see that from $n=50$ simulations, the credibility interval is included in the confidence interval estimated by the naive Monte Carlo method. Furthermore, it decreases slightly with the addition of points to the design of experiments. These results show that adding derivatives to the model allows for significant model improvement and consequently, a more reliable estimate of the probability of failure. In addition, we see in Figure \ref{ggg1} (right) that the uncertainty represented by the credibility intervals is increasing and estimations are less accurate. We have observed this phenomenon several times during the various simulations conducted for the study of this industrial case. This led us to the following hypothesis that a choice of an anisotropic correlation function would be more appropriate. Note that adding degrees of freedom could also benefit to the model including derivative information. However, deriving in this case the covariance matrix requires relatively complex analytical calculations (refer to \cite{lucie_Thesis} for more details). We only performed this burdensome work for the Mat\'ern isotropic correlation function, with $\nu = q + \frac{1}{2}$, $q\in\mathbb{N}^*$. This could be an interesting development perspective.
  


\section{Conclusion and perspectives}\label{Section8}

When a restricted number of data are available, it is common to use a Bayesian approach to estimate the quantity of interest. In our case, this is the product sensitivity to variations in the manufacturing process, expressed as a probability of failure.  Since each observation is costly to obtain, it seems reasonable to build a surrogate probabilistic model (Kriging) and to retain the distribution of the random variable $\Sn$, defined in Equation (\ref{R1}). In this way, we provide an estimation of the failure probability but also quantify the uncertainty introduced by the model. However, this distribution being inaccessible, it is necessary to propose new solutions to enable such an approach. In this article, our objective was therefore to improve the Bayesian estimation procedure by providing practical and efficient means to assess the quality of the prediction. With our main result on the convex order between $\Sn$ and the alternative random variable $\Rn$ given in Proposition \ref{prop_cx_order}, we believe that we achieve this goal. For example, approximations by default and by excess of quantiles of $\Sn$ are deduced. We also show that we can bound from above the moments and the variance of $\Sn$. Moreover, we proved that all these quantities are very easy to estimate using a naive Monte Carlo method. In the same framework, we also derived a sequential procedure to choose a design of experiments based on  the principle of the SUR strategies. We tested our different methods on a real industrial case and proved their interest. Indeed, the conclusion of this analysis is that it is possible to make quantitative risk predictions on real cases, with a very reasonable number of simulations, in a relatively light computer configuration for a company.\medskip

There are many perspectives on this  work. First of all, there are a large number of properties related to the convex order (see \cite{StochasticOrders07}) and others could be used to improve our knowledge of the distribution of $\Sn$ (see, e.g.\cite{Distance_CX_Order}, 
but also \cite{BauerleMuller05} and \cite{MRA} for applications related to risk measures). 
Regarding {SUR} strategies, we proposed a criterion based on the variance of $\Rn$ and, in order to confirm its interest, a full theoretical justification would be interesting, for example, based on the results given in \cite{bect_martingale}. In addition, it would be possible to implement a strategy that selects several points at a time, in order to reduce the number of loops required to converge to the result. The work carried out in \cite{Chevalier14} on this particular subject must be able to be extended. In addition, since {SUR} strategies were not initially developed for vector output, it would be interesting to consider a dedicated approach, as the multi-response case is a common situation in many industrial issues.

\section{Proofs}\label{Section9}	
\subsection{Proof of Proposition \ref{prop_cx_order}}	
We aim at proving that $\Sn\leq_{cx}\Rn$. The proof involves the concept of comonotonicity for a random vector, so we first of all recall the definition of a comonotonic vector (see \cite{Como} and \cite{Kaas} for more information):
\begin{defi}\label{ldodjdmqiq}
A $\mathbb{R}^N$-valued random vector $(Y_i)_{1\leq i \leq N}$ is called comonotonic if, for all $(y_i)_{1\leq i \leq N}\in\mathbb{R}^N$, we have:\begin{align}\label{Def_commonotonic}
\mathbb{P}(Y_1\leq y_1,\ldots, Y_N\leq y_N)=\displaystyle\min_{1\leq i\leq N} \big(\mathbb{P}(Y_i\leq y_i)\big).
\end{align}
\end{defi}
Recall that the convex order definition is given in Definition \ref{def_cx_order}. In the proof, we also use the two following propositions, which can respectively be founded in Section 5 of \cite{Kaas} and Chapter 3 of \cite{StochasticOrders07}: 

\begin{pro}\label{Kaas}
If the random vector $(X_i)_{1\leq i \leq N}$ is comonotonic and has the same marginals as $(Y_i)_{1\leq i \leq N}$, then:\begin{align*} \sum_{i=1}^N Y_i \leq_{cx} \sum_{i=1}^N X_i.\end{align*}
\end{pro}


\begin{pro}\label{CX}
Let $X$ and $Y$ be two random variables, $(X_N)_{N\in\mathbb{N}^*}$ and  $(Y_N)_{N\in\mathbb{N}^*}$ 
be two sequences of random variables such that $X_N\xrightarrow[N\to\infty]{\mathcal{D}}X$ and $Y_N \xrightarrow[N\to\infty]{\mathcal{D}} Y$. 
If the following properties are satisfied: \begin{enumerate}[label=(\roman*), align = right]
\item $\underset{N\to \infty}\lim\mathbb{E} \vert X_N \vert = \mathbb{E}\vert X\vert$ and  $\underset{N\to \infty}\lim\mathbb{E}\vert Y_N\vert = \mathbb{E}\vert Y\vert$, 
\item $X_N \leq_{cx} Y_N$, $\forall N\in\mathbb{N}^*$,
\end{enumerate} then $X \leq_{cx} Y$.
\end{pro}


For all $\x\in\mathbb{X}$, we define the random variable $B_{\x}=\mathbbm{1}_{p_n(\x)>U}$ as a measurable function of the standard uniform random variable $U$. Since $B_{\x}$ has a Bernoulli distribution with parameter $p_n(\x)$, it follows that its cumulative distribution function satisfies:
\begin{align*}
\mathbb{P}(B_\x \leq b) =  \left\{
	\begin{array}{ll} 1 - p_n(\x) & ~\mbox{if}~ b\in[0,1),\\
						 		1 & ~\mbox{if}~ b =  1.
	\end{array}						 
\right. 
\end{align*}

$\rhd$ Application of Definition \ref{ldodjdmqiq}. Let us first prove that for all $N$-tuple $(\x_i)_{1\leq i \leq N}$, the random vector $(B_{\x_{i}})_{1\leq i \leq N}$ is comonotonic, that is: 
\begin{align}\label{Def_commonotonic}
\mathbb{P}(B_{\x_1}\leq b_1,\ldots,B_{\x_N}\leq b_N)=\displaystyle\min_{1\leq i\leq N} \big(\mathbb{P}(B_{\x_i}\leq b_i)\big),
\end{align}
where $b_i\in[0,1]$, for all $1\leq i\leq N$. The random variables $(B_{\x_i})_{1\leq i \leq N}$ taking only the values 0 and 1, it is enough to show it for $b_i\in\lbrace 0,1\rbrace$. For that purpose, we define the finite set $E = \lbrace b_1,\ldots,b_N\rbrace$. Firstly, we consider the particular case where all the $b_i$'s are equal to 1. We then have $ \mathbb{P}(B_{\x_i}\leq b_i)=1$, $\forall i=1,\ldots,N$. This implies: \begin{align*}\mathbb{P}(B_{\x_1}\leq b_1,\ldots,B_{\x_N}\leq b_N)=1=\min_{1\leq i\leq N}\big(\mathbb{P}(B_{\x_i}\leq b_i)\big), 
\end{align*}so that the Equality (\ref{Def_commonotonic}) is satisfied. Secondly, we consider the case where $E$ contains exactly $j$ elements equal to 0, with $j\in\lbrace 1,\ldots,N\rbrace$. By denoting  $\mathfrak{S}(E)$  the symmetric group on $E$, there exists a permutation $\sigma\in\mathfrak{S}(E)$ such that:\\
\begin{center}
$\left\{
\begin{array}{lr}
b_{\sigma(i)}=0, \quad\forall i=1,\ldots,j, \\
b_{\sigma(i)}=1,\quad\forall i=j+1,\ldots,N.
\end{array}\right.$
\end{center}
It follows that:\begin{align*}
 & \mathbb{P}(B_{\x_1}\leq b_1,\ldots,B_{\x_j}\leq b_j, B_{\x_{j+1}}\leq b_{j+1},\ldots, B_{\x_N}\leq b_N)
  \\
 &  = \mathbb{P}(B_{\x_{\sigma(1)}}\leq 0,\ldots,B_{\x_{\sigma(j)}}\leq 0, B_{\x_{\sigma(j+1)}}\leq 1,\ldots, B_{\x_{\sigma(N)}}\leq 1) \\ & = \mathbb{P}(B_{\x_{\sigma(1)}}\leq 0,\ldots,B_{\x_{\sigma(j)}}\leq 0)
 \\ & = \mathbb{P}(p_n(\x_{\sigma(1)})\leq U,\ldots, p_n(\x_{\sigma(j)})\leq U)
 \\ & =1- \displaystyle\max_{1\leq i\leq j}(p_n(\x_{\sigma(i)})),
\end{align*}
because $U$ is uniform on $[0,1]$. Hence, we have:\begin{align*}
\mathbb{P}(B_{\x_1}\leq b_1,\ldots,B_{\x_N}\leq b_N)
= \displaystyle\min_{1\leq i\leq j}\big(1-p_n(\x_{\sigma(i)})\big)
&= \displaystyle\min_{1\leq i\leq j}\big(\mathbb{P}(B_{\x_{\sigma(i)}}\leq 0)\big)= \displaystyle\min_{1\leq i\leq N}\big(\mathbb{P}(B_{\x_{\sigma(i)}}\leq b_{\sigma(i)})\big),
\end{align*}
because $\mathbb{P}(B_{\x_{\sigma(i)}}\leq b_{\sigma(i)})=1$,  $\forall i=j+1,\ldots,N.$ Finally, we have:\begin{align*}
\mathbb{P}(B_{\x_1}\leq b_1,\ldots,B_{\x_N}\leq b_N) & =\displaystyle\min_{1\leq i\leq N}\big(\mathbb{P} (B_{\x_{\sigma(i)}}\leq b_{\sigma(i)})\big)  =\displaystyle\min_{1\leq i\leq N}\big(\mathbb{P}(B_{\x_i}\leq b_i)\big),
\end{align*}so that Equality (\ref{Def_commonotonic}) is verified again. We then proved that, for all $(\x_i)_{1\leq i \leq N}\in\mathbb{X}^N$, the vector $(B_{\x_{i}})_{1\leq i \leq N}$ is comonotonic.\medskip

$\rhd$ \textit{Application of Theorem \ref{Kaas}.} Since, for all $n$-tuple $(\x_1,\ldots,\x_N)\in\mathbb{X}^N$, the vector $(B_{\x_{i}})_{1\leq i \leq N}=(\mathbbm{1}_{p_n(\x_i)>U})_{1\leq i \leq N}$ is comonotonic and has the same marginal distributions as the vector $ (\mathbbm{1}{_{\xi_n(\x_i)> T}})_{1\leq i \leq N}$, we can apply Theorem \ref{Kaas}. As a result, we have:\begin{align*}
\sum_{i=1}^N\mathbbm{1}{_{\xi_n(\x_i)> T}}\leq_{cx} \sum_{i=1}^N \mathbbm{1}_{p_n(\x_i)>U},
\end{align*}
or equivalently,
\begin{align}\label{sum_cx_order}
\dfrac{1}{N}\sum_{i=1}^N\mathbbm{1}{_{\xi_n(\x_i)> T}}\leq_{cx} \dfrac{1}{N}\sum_{i=1}^N \mathbbm{1}_{p_n(\x_i)>U}.
\end{align}

$\rhd$ \textit{Application of Proposition \ref{CX}}. Let us consider a random sample $\xx_1,\ldots,\xx_N$ of i.i.d. copies of $\xx$ distributed with respect to $\loiX$. We define:
$$V_N=\frac{1}{N}\displaystyle\sum_{i=1}^N\mathbbm{1}_{\xi_n(\xx_i)> T}\quad \mbox{and}\quad W_N=\frac{1}{N}\displaystyle\sum_{i=1}^N\mathbbm{1}{_{p_n(\xx_i)>U}}.$$ 
According to the law of large numbers, we have:\begin{align*} 
V_N\xrightarrow[N\to\infty]{a.s.}
\mathbb{E}[\mathbbm{1}_{\xi_n(\xx)>T}\mid\xi_n] = \Sn
\quad \mbox{and}\quad 
W_N\xrightarrow[N\to\infty]{a.s.}\mathbb{E}[\mathbbm{1}_{p_n(\xx)>U}\mid U]=\Rn.
\end{align*} 
Let $h$ be a continuous function. Then by the continuity theorem, knowing $\xi_n$, we have:\begin{align*} 
h(V_N)\xrightarrow[N\to\infty]{a.s.} h(\Sn),
\end{align*}and, knowing $U$, we have:\begin{align*} 
h(W_N)\xrightarrow[N\to\infty]{a.s.} h(\Rn).
\end{align*}
Since all the above variables take  values in $[0,1]$, then the expectations exist and by Lebesgue's (conditional) dominated convergence theorem, we have:\begin{align*} 
\mathbb{E}[h(V_N)\mid \xi_n]\xrightarrow[N\to\infty]{a.s.}\mathbb{E}[ h(\Sn)\mid \xi_n]\quad \mbox{and}\quad 
\mathbb{E}[h(W_N)\mid U]\xrightarrow[N\to\infty]{a.s.}\mathbb{E}[ h(\Rn)\mid U].
\end{align*}This implies that:\begin{align}\label{oi} 
\lim_{N\to \infty}\mathbb{E}[h(V_N)]=\mathbb{E}[h(\Sn)]
\quad \mbox{and}\quad 
\lim_{N\to \infty}\mathbb{E}[h(W_N)]=\mathbb{E}[h(\Rn)].
\end{align}As a result, we have: \begin{align*} V_N\xrightarrow[N\to\infty]{\mathcal{L}}\Sn\quad\mbox{and}\quad W_N\xrightarrow[N\to\infty]{\mathcal{L}}\Rn.\end{align*}In particular, by taking $h(x) = \vert x\vert$ in (\ref{oi}), the assumption \textit{(i)} of Proposition \ref{CX} is satisfied. Moreover, by simultaneously taking into account the Definition \ref{def_cx_order} of a convex order and Inequality (\ref{sum_cx_order}), for any convex function $\varphi$, we have: \begin{align*}\mathbb{E}\left[\varphi\left(V_N\right)\mid \xx_1,\ldots,\xx_N \right]\leq \mathbb{E}\left[\varphi\left(W_N\right)\mid \xx_1,\ldots,\xx_N \right].\end{align*}
This implies that:\begin{align*}\mathbb{E}\left[\varphi(V_N)\right]\leq\mathbb{E}\left[\varphi(W_N)\right],\end{align*}
or, equivalently, \begin{align*}V_N\leq_{cx}W_N.\end{align*}Thus, the hypothesis \textit{(ii)} of Proposition \ref{CX} is also verified and this proposition applies. It follows that $\Sn\leq_{cx}\Rn$.

	\subsection*{Proof of Proposition \ref{Prop_CDF_Quantile_function_Rn}} \label{Demo_Prop_CDF_Quantile_function_Rn}
	
We aim at proving that for all $\alpha\in[0,1]$, the $\alpha$-quantile $F_{\Rn}^{-1}(\alpha)$ satisfies:\begin{align*}
F_{\Rn}^{-1}(\alpha) = \int_{\set}\mathbbm{1}_{p_n(\x)>1-\alpha}\loiX(d\x).
\end{align*} First of all, let us introduce some notations. Let $F_{\Rn}$ be the cumulative distribution function of $\Rn$ and $G_{\Rn} = 1- F_{\Rn}$ be its survival function. We recall that $F_{\Rn}^{-1}(\alpha) $ is defined by:\begin{align*}
F_{\Rn}^{-1}(\alpha)  = \inf\lbrace t\in[0,1] : F_{\Rn}(t) \geq \alpha\rbrace
&  = \inf\lbrace t\in[0,1] : G_{\Rn}(t) \leq 1-\alpha\rbrace.
\end{align*}Moreover, let $G_{p_n(\xx)}$ be the survival function of the random variable $p_n(\xx)$:\begin{align*}
G_{p_n(\xx)}(u) = \mathbb{P}(p_n(\xx)>u) = \int_{\set}\mathbbm{1}_{p_n(\x)>u}\loiX(d\x), \quad\forall u \in [0,1].
\end{align*}By taking $U$ as a random variable with uniform distribution on $[0,1]$, it immediately comes that: \begin{align*} G_{p_n(\xx)}(U)= \mathbb{P}(p_n(\xx)>U\mid U) =\Rn.\end{align*} Now, let us consider the function $G_{p_n(\xx)}^{-1}$ defined for all $t\in[0,1]$ by: 
 \begin{align*}
G_{p_n(\xx)}^{-1}(t) = \inf\lbrace u\in[0,1] : G_{p_n(\xx)}(u) \leq t\rbrace.\end{align*}It can be verified that functions $G_{p_n(\xx)}$ and $G_{p_n(\xx)}^{-1}$ satisfy:\begin{align*}G_{p_n(\xx)}(u) \leq t \Leftrightarrow u\geq G_{p_n(\xx)}^{-1}(t), \quad\forall u\in[0,1]\mbox{ and }\forall t \in [0,1].\end{align*}
For all $t\in[0,1]$, the survival function $G_{\Rn}$ then satisfies:  \begin{align*}
G_{\Rn}(t) = \mathbb{P}(\Rn> t) = \mathbb{P}(G_{p_n(\xx)}(U)> t) = \mathbb{P}(U < G_{p_n(\xx)}^{-1}(t)) =  G_{p_n(\xx)}^{-1}(t).
\end{align*}Finally, for all $\alpha\in [0,1]$, we have:\begin{align*}
F_{\Rn}^{-1}(\alpha) = \inf\lbrace t\in[0,1] : G_{\Rn}(t) \leq 1-\alpha\rbrace  & = \inf\lbrace t\in[0,1] : G_{p_n(\xx)}^{-1}(t) \leq 1-\alpha\rbrace\\
&  = \inf\lbrace t\in[0,1] : G_{p_n(\xx)}(1-\alpha) \leq t\rbrace\\
& = G_{p_n(\xx)}(1-\alpha)\\
& = \int_\set\mathbbm{1}_{p_n(\x)>1-\alpha}\loiX(d\x).
\end{align*}

	\subsection*{Proof of Proposition \ref{Prop_CredibleInt_Markov}} \label{Demo_Prop_CredibleInt_Markov}
	
Recall that $\mu_n = \mathbb{E}[\Rn] = \displaystyle\int_0^1 F_{\Rn}^{-1}(t)dt$. Here, we aim at proving that for all $\alpha\in(0,1)$, we have:\begin{align}\label{iiiii}
\frac{\mu_n + \alpha-1}{\alpha}\leq\dfrac{1}{\alpha}\int_0^\alpha F_{\Rn}^{-1}(t)dt\leq F_{\Sn}^{-1}(\alpha) \leq\dfrac{1}{1-\alpha}\int_{\alpha}^1F_{\Rn}^{-1}(t)dt \leq \dfrac{\mu_n}{1-\alpha},
\end{align}According to Proposition \ref{prop_cx_order}, we have $\Sn\leq_{cx}\Rn$. Thus, for all $\alpha\in(0,1)$, we have (see Theorem 3.A.5. in \cite{StochasticOrders07}):
\begin{align*}
\int_{0}^\alpha F_{\Sn}^{-1}(t)\emph{\mbox{d}}t\geq\int_{0}^\alpha F_{\Rn}^{-1}(t)\emph{\mbox{d}}t\quad\mbox{and}\quad
\int_{\alpha}^1F_{\Sn}^{-1}(t)\emph{\mbox{d}}t\leq\int_{\alpha}^1F_{\Rn}^{-1}(t)\emph{\mbox{d}}t.
\end{align*}
The function $F_{\Sn}^{-1}$ being monotonically increasing, we deduce that: \begin{align*}
F_{\Sn}^{-1}(\alpha) & = \dfrac{1}{\alpha}\displaystyle\int_0^\alpha F_{\Sn}^{-1}(\alpha)dt  \geq \dfrac{1}{\alpha}\displaystyle\int_0^\alpha F_{\Sn}^{-1}(t)dt  \geq \dfrac{1}{\alpha}\displaystyle\int_0^\alpha F_{\Rn}^{-1}(t)dt.\end{align*}
Yet, $\mu_n  = \displaystyle\int_0^1 F_{\Rn}^{-1}(t)dt  \leq  \displaystyle\int_0^\alpha F_{\Rn}^{-1}(t)dt + 1-\alpha$,
because $0\leq F_{\Rn}^{-1}(t)\leq 1$, $\forall t\in [0,1]$. Consequently, \begin{align*}\mu_n + \alpha -1 \leq \displaystyle\int_0^\alpha F_{\Rn}^{-1}(t)dt\quad\mbox{and}\quad\frac{\mu_n + \alpha -1}{\alpha}  \leq \dfrac{1}{\alpha}\displaystyle\int_0^\alpha F_{\Rn}^{-1}(t)dt\leq F_{\Sn}^{-1}(\alpha).
\end{align*}
Similarly, we have:\begin{align*}
F_{\Sn}^{-1}(\alpha)  = \dfrac{1}{1 - \alpha}\displaystyle\int_\alpha^1 F_{\Sn}^{-1}(\alpha)dt  \leq \dfrac{1}{1 - \alpha}\displaystyle\int_\alpha^1 F_{\Sn}^{-1}(t)dt \leq \dfrac{1}{1 - \alpha}\displaystyle\int_\alpha^1 F_{\Rn}^{-1}(t)dt,
\end{align*}
and \begin{align*}\dfrac{1}{1 - \alpha}\displaystyle\int_\alpha^1 F_{\Rn}^{-1}(t)dt\leq \dfrac{1}{1-\alpha}\displaystyle\int_0^1 F_{\Rn}^{-1}(t)dt = \dfrac{\mu_n}{1-\alpha},\end{align*}
 
		\subsection*{Proof of Proposition \ref{Exp_bounds}} \label{Demo_Exp_bounds}

According to Proposition \ref{Prop_CDF_Quantile_function_Rn}, we have $
F^{-1}_{\Rn}(\alpha) = \displaystyle\int_{\mathbb{X}}\mathbbm{1}_{\alpha>1 - p_n(\x)}\textbf{P}_{\xx}(d\x)$, $\forall\alpha\in [0,1]$.
As a result:
\begin{align*}
\dfrac{1}{\alpha}\int_0^\alpha F^{-1}_{\Rn}(t)dt & = \dfrac{1}{\alpha}\displaystyle\int_{\mathbb{X}}\left(\int_0^\alpha\mathbbm{1}_{t > 1 - p_n(\x)}dt\right)\textbf{P}_{\xx}(d\x)\\
& = \dfrac{1}{\alpha}\displaystyle\int_{\mathbb{X}}\max\big(0, \alpha - 1 + p_n(\x)\big)\textbf{P}_{\xx}(d\x)\\
&  = 1 -  \displaystyle\int_{\mathbb{X}}\min\left(1, \dfrac{1 - p_n(\x)}{\alpha}\right)\textbf{P}_{\xx}(d\x), 
\end{align*}
and
\begin{align*}
\dfrac{1}{1 - \alpha}\int_\alpha^1 F^{-1}_{\Rn}(t)dt & = \dfrac{1}{1 - \alpha}\displaystyle\int_{\mathbb{X}}\left(\int_\alpha^1\mathbbm{1}_{t > 1 - p_n(\x)}dt\right)\textbf{P}_{\xx}(d\x)\\
& = \dfrac{1}{1 - \alpha}\displaystyle\int_{\mathbb{X}}\min\big(1 - \alpha , p_n(\x)\big)\textbf{P}_{\xx}(d\x)\\
&  = \displaystyle\int_{\mathbb{X}}\min\left(1,\dfrac{p_n(\x)}{1 - \alpha}\right)\textbf{P}_{\xx}(d\x). 
\end{align*}

\subsection*{Proof of Proposition \ref{Inequality_Criteria}} \label{SUR_inegality_criterion}

The following demonstration is based on the one for Proposition 3 in \cite{Li12}, which shows that $\JS\leq J_{n,k}$, $\forall k = 1, \ldots, 4$. In the following, we specify this partial order  and show that $\JR$ offers a better local approximation of $\JS$ than criteria $( J_{n,k})_{k = 1, \ldots, 4}$.\medskip

According to the convex order inequality established in the Proposition \ref{prop_cx_order}, we have $\mbox{Var}[\Sn]\leq\mbox{Var}[\Rn]$. This means that $\JS(\x)\leq\JR(\x)$, $\forall\x\in\mathbb{X}$. Moreover, by definition of $\Rn$, we have:
\begin{equation*}
\Rn-\mathbb{E}[\Rn]=\displaystyle\int_{\mathbb{X}}\left(\mathbbm{1}_{p_n(\x)> U}-p_n(\x)\right)\textbf{P}_{\xx}(d\x). 
\end{equation*}
Let us denote by $\vert X\vert = {\mathbb{E}[X^2]}^\frac{1}{2}$ the Euclidean norm defined on the space $L^2(\Omega, \mathcal{F}, \mathbb{P})$ of the random integrable square variables. For all $\x\in\E$, based on the generalized Minkowski inequality (see \cite{Minkowski}) and given that $
\int_0^1\mathbbm{1}_{p_n(\x)>u}du = p_n(\x)$, we have:\begin{align*}\Vert\Rn-\mathbb{E}[\Rn]\Vert & = \left(\displaystyle\int _0^1\left(\displaystyle\int_{\E}\mathbbm{1}_{p_n(\x)>u}-p_n(\x)\textbf{P}_{\xx}(d\x)\right)^2{d}u\right)^{\frac{1}{2}}\\
								& \leq \int_{\E}\bigg(\int_0^1\big( \mathbbm{1}_{p_n(\x)> u}-p_n(\x)\big)^2{d}u\bigg)^{\frac{1}{2}}\textbf{P}_{\xx}(d\x)\\
                        		 & =\int_{\E}\big(p_n(\x)(1-p_n(\x)\big)^\frac{1}{2}\textbf{P}_{\xx}(d\x).
\end{align*}Then,
\begin{align*}
\mbox{Var}[\Rn]  =\Vert \Rn-\mathbb{E}[\Rn]\Vert^2 & \leq\left(\int_{\E}\big(p_n(\x)(1-p_n(\x)\big)^\frac{1}{2}\textbf{P}_{\xx}(d\x)\right)^2  \leq\int_{\E}p_n(\x)(1-p_n(\x))\textbf{P}_{\xx}(d\x),
\end{align*} according to Jensen's inequality. 
In addition, since for all $x\in [0,1]$, we have $x(1-x)\leq \min(x, 1-x)$, it follows that:\begin{align*}
\mbox{Var}[\Rn]  \leq \int_{\E} \min(p_n(\x), 1-p_n(\x))\textbf{P}_{\xx}(d\x).
\end{align*}We then proved that:
\begin{align*}\JR(\x)\leq J_{n,2}(\x)\leq J_{n,4}(\x)\leq J_{n,3}(\x), \quad\forall\x\in\mathbb{X}.\end{align*}
We also have:
\begin{align*}
\mbox{Var}[\Rn] & \leq\left(\int_{\E}\big(p_n(\x)(1-p_n(\x)\big)^\frac{1}{2}\textbf{P}_{\xx}(d\x)\right)^2 \leq\left(\int_{\E}\min\big(p_n(\x, 1-p_n(\x)\big)^\frac{1}{2}\textbf{P}_{\xx}(d\x)\right)^2, 
\end{align*}that is:\begin{align*}\JR(\x)\leq J_{n,2}(\x)\leq J_{n,1}(\x), \quad\forall\x\in\mathbb{X}.\end{align*}	We then proved the following partial order:\begin{align*}
\JS\leq\JR\leq J_{n, k},\quad \forall k=1, \ldots,4.
\end{align*}

	\subsection*{Proof of Proposition \ref{varR2_ClassError}}\label{qrfAAAA} 

By definition of $\Rn$, its variance satisfies:
\begin{align*}
\mbox{Var}[\Rn] &  = \mathbb{E}\bigg[\bigg(\int_{\E}\mathbbm{1}_{p_n(\x)>U}\loiX(d\x) - \int_{\E}p_n(\x)\loiX(d\x)\bigg)^2\bigg] \\
& =\int_{\E^2}\mathbb{E}\big[\mathbbm{1}_{p_n(\x)>U}\mathbbm{1}_{p_n(\y)>U}\big]\loiX(d\x)\loiX(d\y) - \int_{\E^2}p_n(\x)p_n(\y)\loiX(d\x)\loiX(d\y)  \\
& = \displaystyle\int_{\mathbb{X}^2}\big(\min(p_n(\x), p_n(\y))-p_n(\x)p_n(\y)\big)\textbf{P}_{\xx}(d\x)\textbf{P}_{\xx}(d\y),\\
& = \displaystyle\int_{\mathbb{X}^2}\left((1-p_n(\x))p_n(\y)\mathbbm{1}_{p_n(\x)\geq p_n(\y)}+p_n(\x)(1-p_n(\y))\mathbbm{1}_{p_n(\x)< p_n(\y)}\right)\textbf{P}_{\xx}(d\x)\textbf{P}_{\xx}(d\y),\\
& =  \displaystyle\int_{\mathbb{X}}\bigg((1-p_n(\x))\displaystyle\int_\E p_n(\y)\mathbbm{1}_{p_n(\x)\geq p_n(\y)}\textbf{P}_{\xx}(d\y)+p_n(\x)\displaystyle\int_{\E}(1-p_n(\y))\mathbbm{1}_{p_n(\x)< p_n(\y)}\textbf{P}_{\xx}(d\y)\bigg)\textbf{P}_{\xx}(d\x)\\
& = \int_\set \eta_n(\x)\loiX(d\x),
\end{align*}where $\eta_n(\x) = (1-p_n(\x))\displaystyle\int_\E p_n(\y)\mathbbm{1}_{p_n(\x)\geq p_n(\y)}\textbf{P}_{\xx}(d\y) + p_n(\x)\displaystyle\int_{\E}(1-p_n(\y))\mathbbm{1}_{p_n(\x)< p_n(\y)}\textbf{P}_{\xx}(d\y)$. 
	
	\subsection*{Proof of Proposition \ref{Max_etaFunc}}\label{proffffff} 
	
In this demonstration, we assume that the random variable $p_n(\xx)$ is absolutely continuous, i.e. it admits a density $f_{p_n(\xx)}$.  We denote by $G_{p_n(\xx)}$ its survival function and $F_{p_n(\xx)}$ its cumulative distribution function:\begin{align*}G_{p_n(\xx)}(t) = \mathbb{P}(p_n(\xx)>t) = \int_0^1\mathbbm{1}_{u>t}f_{p_n(\xx)}(u)du = \int_\E\mathbbm{1}_{p_n(\x)>t}\loiX(d\x) = 1 - F_{p_n(\xx)}(t),\quad\forall t\in [0,1].\end{align*}
The expectation $\mu_n$ of $\Sn$ is then written: \begin{align*}
\mu_n = \mathbb{E}[\Sn] = \mathbb{E}[p_n(\xx)] = \int_{0}^1 u f_{p_n(\xx)}(u)du = \int_0^1\mathbb{P}(p_n(\xx)>u)du = \int_0^1G_{p_n(\xx)}(u)du.
\end{align*}
Let us consider the function $\eta_n$ given in Proposition \ref{varR2_ClassError} and defined for all $\x\in\mathbb{X}$ by: \begin{align}\eta_n(\x) = p_n(\x)\displaystyle\int_{\E}\big(1-p_n(\y)\big)\mathbbm{1}_{p_n(\x)< p_n(\y)}\textbf{P}_{\xx}(d\y) + \big(1-p_n(\x)\big)\displaystyle\int_\E p_n(\y)\mathbbm{1}_{p_n(\x)\geq p_n(\y)}\textbf{P}_{\xx}(d\y).
\end{align}We will show that this function has a global maximum on $[0,1]$ at the point: $$q_n^*=\mathbb{P}(\Rn>\mu_n).$$ To do this, let us start by noting that $\eta_n$ can be rewritten as follow:\begin{align*}\eta_n(\x)  & = p_n(\x)\mathbb{E}\big[(1-p_n(\xx))\mathbbm{1}_{p_n(\x)\leq p_n(\xx)}\big] + \big(1-p_n(\x)\big)\mathbb{E}\big[p_n(\xx)\mathbbm{1}_{p_n(\x)\geq p_n(\xx)}\big] \\ 
& =  p_n(\x)\displaystyle\int_{p_n(\x)}^1(1-u)f_{p_n(\xx)}(u)du + \big(1-p_n(\x)\big)\int_0^{p_n(\x)}uf_{p_n(\xx)}(u)du\\
& = p_n(\x)\int_{p_n(\x)}^1 f_{p_n(\xx)}(u)du - \mu_ np_n(\x) + \int_0^{p_n(\x)} u f_{p_n(\xx)}(u)du\end{align*}
By applying an integration by parts to the right term, we obtain:
\begin{align*}
\eta_n(\x) & = p_n(\x)\big(1 - F_{p_n(\xx)}(p_n(\x))\big) - \mu_ np_n(\x) +  \Big[p_n(\x)F_{p_n(\xx)}(p_n(\x)) - \int_0^{p_n(\x)}F_{p_n(\xx)}(u)du\Big]\\
&   = \int_{0}^{p_n(\x)}G_{p_n(\xx)}(u)du - \mu_ np_n(\x).
\end{align*}	
We set $\varphi(q) =  \displaystyle\int_{0}^{q}G_{p_n(\xx)}(u)du - \mu_n q$, $\forall q\in [0,1]$. Then,
\begin{align}\label{cvcgd}\varphi'(q) = 0 \Leftrightarrow G_{p_n}(q) = \mu_n\Leftrightarrow\int_{\mathbb{X}}\mathbbm{1}_{p_n(\x)>q}\textbf{P}_{\xx}(d\x) = \mu_n.	
\end{align} 			
Let $G_{p_n(\xx)}^{-1}$ be the function defined for all $t\in[0,1]$ by $G_{p_n(\xx)}^{-1}(t) = \inf\lbrace u \in[0,1] : G_{p_n(\xx)}(u)\leq t\rbrace$.
As mentioned in proof of Proposition \ref{Prop_CDF_Quantile_function_Rn}, we have: 
\begin{align*}G_{p_n(\xx)}^{-1}(t) = \mathbb{P}(\Rn>t),\quad\mbox{and}\quad
G_{p_n(\xx)}(u)\leq t \Leftrightarrow u \geq G_{p_n(\xx)}^{-1}(t), \quad\forall u\in[0,1],~\forall t\in[0,1]. 
\end{align*}Then, $
G_{p_n(\xx)}(q) = \mu_n\Leftrightarrow q = \mathbb{P}(\Rn>\mu_n)$ and $\eta_n$ has a unique extrema at point $q_n^* = \mathbb{P}(\Rn>\mu_n)$. According to (\ref{cvcgd}), it satisfies: $\int_{\mathbb{X}}\mathbbm{1}_{p_n(\x)>q_n^*}\textbf{P}_{\xx}(d\x) = \mu_n$. By studying the sign of $\varphi'$, it is easy to verify that it corresponds to a maximum.

\section*{Acknowledgments}

The authors thank Arnaud Guyader and Florent Malrieu for their help on several aspect of this paper and for proofreading.

\bibliographystyle{abbrv} 
\bibliography{BiblioManuscrit}

\end{document}